\documentclass[12pt,a4paper]{article}
\usepackage[english]{babel}
\usepackage{ae,aecompl}
\usepackage{amssymb,amsthm,amsmath}
\usepackage{fullpage}
\title{The noncommutative geometry \\[5pt] of the quantum projective plane}

\date{20 December 2007; revised 8 September 2008}

\author{~\\
\large{Francesco D'Andrea$^1$, Ludwik D\k{a}browski$^2$, Giovanni Landi$^3$} \\ [20pt]
 \normalsize{$^1$ D{\'e}partement de Math{\'e}matique, Universit\'e Catholique de Louvain,} \\
 \normalsize{Chemin du Cyclotron 2, B-1348, Louvain-La-Neuve, Belgium} \\ [10pt]
 \normalsize{$^2$ Scuola Internazionale Superiore di Studi Avanzati,}\\
 \normalsize{Via Beirut 2-4, I-34014, Trieste, Italy} \\ [10pt]
 \normalsize{$^3$ Dipartimento di Matematica e Informatica, Universit{\`a} di Trieste,} \\
\normalsize{Via A.~Valerio 12/1, I-34127, Trieste, Italy} \\
\normalsize{and INFN, Sezione di Trieste, Trieste, Italy}
\\~\\
}

\allowdisplaybreaks[4]
\newenvironment{prova}{\begin{proof}\parindent=0in}{\end{proof}}
\numberwithin{equation}{section}

\newtheorem{thm}{Theorem}[section]
\newtheorem{lemma}[thm]{Lemma}
\newtheorem{prop}[thm]{Proposition}

\usepackage{hyperref}
\hypersetup{%
 colorlinks=true,
 pdfstartpage=1,
 pdfstartview=FitH,
 bookmarksnumbered=true,
 linkcolor=black,
 anchorcolor=black,
 citecolor=black,
 urlcolor=black}

\newcommand{\arxiv}[1]{[arxiv:\htmladdnormallink{#1}{http://arxiv.org/abs/#1}]}
\newcommand{\doi}[1]{[doi:\htmladdnormallink{#1}{http://dx.doi.org/#1}]}

\newcommand{\A}{\mathcal{A}}
\newcommand{\B}{\mathcal{B}}
\newcommand{\U}{\mathcal{U}}
\newcommand{\HH}{\mathcal{H}}
\newcommand{\N}{\mathbb{N}}
\newcommand{\Z}{\mathbb{Z}}
\newcommand{\R}{\mathbb{R}}
\newcommand{\C}{\mathbb{C}}
\newcommand{\CP}{\mathbb{C}\mathrm{P}}
\newcommand{\oh}{\smash[t]{\smash[b]{\tfrac{1}{2}}}}
\newcommand{\ma}[1]{\left(\!\begin{array}{cccc}#1\end{array}\!\right)}
\newcommand{\maa}[1]{\text{\footnotesize $\left(\!\begin{array}{ccc}#1\end{array}\!\right)$}}

\newcommand{\sqbn}[1]{\textrm{\footnotesize$\left[\!\!\begin{array}{c}#1\end{array}\!\!\right]$}}
\newcommand{\az}{\triangleright}
\newcommand{\za}{\triangleleft}
\newcommand{\aaz}{\,\textrm{\footnotesize$\blacktriangleright$}\,}
\newcommand{\inner}[1]{\left<#1\right>}
\newcommand{\ket}[1]{\left|#1\right>}

\newcommand{\tr}{\mathrm{Tr}}
\newcommand{\de}{\mathrm{d}}
\newcommand{\wprod}{\wedge_q\!}
\newcommand{\deb}{\bar{\partial}}
\newcommand{\D}{D}
\newcommand\ii{\mathrm{i}}

\begin{document}

\maketitle

\bigskip

\thispagestyle{empty}
                                                                               
\begin{abstract}\noindent
 We  study the spectral geometry of the quantum projective plane $\CP^2_q$, 
 a deformation of the complex projective plane $\CP^2$, 
 the simplest example of a $\mathrm{spin}^c$ manifold which is not spin. 
In particular, we construct a Dirac operator $D$ which gives a $0^+$-summable spectral triple, 
equivariant under $U_q(su(3))$. The square of $D$ is a central 
element for which left and right actions on spinors coincide, a fact 
that is exploited to compute explicitly its spectrum.
\end{abstract}

\vfill



\section{Introduction}

The geometry of quantum spaces
-- whose coordinate algebras are noncommutative --
can be studied, following A. Connes \cite{Con94}, by means of a spectral triple.
The latter is the datum $(\A,\HH,D)$, where $\A$ is a unital, involutive,
associative (but non necessarily commutative) $\C$-algebra 
with a faithful representation, $\pi:\A\to\B(\HH)$, on a separable Hilbert space $\HH$,
and $D$ is a selfadjoint operator on $\HH$ with 
compact resolvent and such that $[D,a]$ is bounded for all $a\in\A$. The operator $D$ is called (a generalized) Dirac operator. In addition, the spectral triple is called \emph{even} if $\HH=\HH_+\oplus\HH_-$ is $\Z_2$-graded,
the representation of $\A$ is diagonal and the operator $D$ is off-diagonal for this decomposition. The requirement of compact resolvent for the Dirac operator guarantees, for example, that in the even case the twisting of the operator $D^\pm=D|_{\HH_\pm}$ with projections (describing classes in the K-theory of $\A$) are unbounded Fredholm operators: the starting
point for the construction of `topological' invariants via index computations.
Roughly, the bounded commutators condition says that the specrum of $D$ does not grow too rapidly, while the compact resolvent one says that the specrum of $D$ does not  grow too slowly. It is the interplay of the two that (together with further requirements) imposes stringent
restrictions on the geometry and produces spectacular consequences.

For quantum homogeneous spaces (that is spaces which are `homogeneous' for quantum groups, see e.g.~\cite{KS97}), 
a possible strategy consists to define a Dirac operator by its spectrum, 
in a suitable basis of `harmonic' spinors, and to prove that the commutators $[D,a]$ are bounded by the use 
of quantum groups representation theory. 
In this manner one usually finds Dirac operators with spectrum growing at most
polynomially (cf.~\cite{CP03a,DLS05,DLPS05,DDLW07}). 
 
A different occurrence is for the standard Podle\'s quantum sphere where also a Dirac operator exists \cite{DS03} with a spectrum growing exponentially, defining then a $0^+$-summable spectral triple (a behaviour on the opposite hand to that of theta-summability).
This operator has a particular geometrical meaning as it can be constructed \cite{SW04} by using the action of certain generators 
of $U_q(su(2))$
which act as derivations on the standard Podle\'s sphere. 
Along this line, a general construction of Dirac operators $D$ on
quantum irreducible flag manifolds, including 
projective spaces, was given in \cite{Kra04}. 
These operators were used to realize the differential calculi of
\cite{HK04} by expressing the exterior derivative as a commutator with $D$.
However, in \cite{Kra04} there is no computation of any spectrum of $D$ and thus no addressing, among other things, of the compact resolvent requirement for the Dirac operator, an essential feature of spectral triples as mentioned above. Furthermore, the construction there depends on the choice of a morphism $\gamma$ (Prop.~2 in \cite{Kra04}) that appears to be neither unique nor canonical. 

In the present paper, as a first step for a general strategy, we work out from scratch the spectral geometry of a basic example (besides the standard Podle\'s sphere), that is the quantum complex projective plane $\CP^2_q$. This is defined as a $q$-deformation with real parameter (that we restrict to $q\in(0,1)$) of the complex projective plane $\CP^2$ seen as the four dimensional real manifold 
$S^5/S^1=SU(3) / SU(2)\times U(1)$. Our example is particularly important in that it is a deformation of a manifold which is not 
a spin manifold but only $\mathrm{spin}^c$. In analogy with the standard Podle\'s sphere, we find a Dirac operator $D$ on $\CP^2_q$ with exponentially growing spectrum -- a $q$-deformation of the spectrum of the Dolbeault-Dirac operator on undeformed $\CP^2$ (for the latter cf. \cite{GS99}) --, thus giving a $0^+$-dimensional spectral triple. The spectrum is explicitly computed by relating the square of $D$ to a 
quantum Casimir element whose left and right actions on spinors coincide. As motivated in Sect.~\ref{sec:2}, to get this quantum Casimir element we need to enlarge the symmetry algebra. 
The use of this technique to compute the spectrum via left/right actions seem to be, to the best of our knowledge, a novel one. 
There remains open problems, notably the issue of regularity for the present spectral geometry, which might hold at most in the `twisted sense' of \cite{NT}; 
their analysis is postponed to future work.

The plan of the paper is the following. In Sect.~\ref{sec:2} we introduce
the Hopf algebra $U_q(su(3))$, 
which describes the `infinitesimal' symmetries of $\CP^2_q$, and in Sect.~\ref{sec:3}
the dual Hopf algebra $\A(SU_q(3))$, whose elements are representative
functions on the quantum $SU(3)$ group. The coordinate algebra of $\CP^2_q$ is
defined in Sect.~\ref{sec:4} as the fixed point
subalgebra of $\A(SU_q(3))$ for the action of a suitable Hopf subalgebra $U_q(u(2))
\subset U_q(su(3))$.
In Sect.~\ref{sec:5} we describe the $q$-analogue of antiholomorphic forms
and use them to construct first a differential calculus and then a spectral triple
on $\CP^2_q$ in Sect.~\ref{sec:spinc}. The appendix contains the description of antiholomorphic forms on the classical $\CP^2$ as equivariant maps, a description which was the motivation for an analogous identification on the quantum $\CP^2_q$. 

\section{The symmetry Hopf algebra $U_q(su(3))$}\label{sec:2}
Let  $U_q(su(3))$ be the Hopf
$*$-algebra generated (as a $*$-algebra) by $ K_i , K_i^{-1},E_i, F_i  \, , \, i=1,2$,
with $ K_i=K_i^*,F_i=E_i^*$, and relations
\begin{gather*}
[K_i,K_j]=0\;,\qquad
K_iE_iK_i^{-1}=qE_i\;,\qquad
[E_i,F_i]=(q-q^{-1})^{-1}(K_i^2-K_i^{-2})\, \\
\rule{0pt}{14pt}
K_iE_jK_i^{-1}=q^{-1/2}E_j\;,\qquad
[E_i,F_j]=0\;, \qquad\mathrm{if}\;i\neq j\;,
\end{gather*}
and (Serre relations)
\begin{equation}\label{serre}
E_i^2E_j-(q+q^{-1})E_iE_jE_i+E_jE_i^2=0\qquad\forall\;i\neq j\;.
\end{equation}
We can restrict the real deformation parameter to the interval $0<q<1$; for $q> 1$ we get isomorphic algebras. In App.~\ref{sec:A} we shall also briefly decribe the `classical limit' $U(su(3))$. In the notation of~\cite[Sect.~6.1.2]{KS97} the above Hopf algebra is denoted $\breve{U}_q(su(3))$, the `compact' real form of the Hopf algebra denoted $\breve{U}_q(sl(3,\C))$ there.
With the $q$-commutator defined as
$$
[a,b]_q:=ab-q^{-1}ba \;,
$$
 relations \eqref{serre} can be rewritten as $\,[E_i,[E_j,E_i]_q]_q=0\,$
or $\,[[E_i,E_j]_q,E_i]_q=0\,$.
\\
Coproduct, counit and antipode are given by (with $i=1,2$)
\begin{gather*}
\Delta(K_i)=K_i\otimes K_i\;,\qquad
\Delta(E_i)=E_i\otimes K_i+K_i^{-1}\otimes E_i\;, \\
\epsilon(K_i)=1\;,\qquad
\epsilon(E_i)=0\;,\qquad
S(K_i)=K_i^{-1}\;,\qquad
S(E_i)=-qE_i\;.
\end{gather*}
The opposite Hopf $*$-algebra $U_q(su(3))^{\mathrm{op}}$ is defined to be
isomorphic to $U_q(su(3))$ as $*$-coalgebra, but equipped with opposite
multiplication and with antipode $S^{-1}$. There is a
Hopf $*$-algebra isomorphism $\vartheta:U_q(su(3))\to U_q(su(3))^{\mathrm{op}}$
given on generators by
\begin{equation}\label{eq:vartheta}
\vartheta(K_i):=K_i\;,\qquad
\vartheta(E_i):=F_i\;,\qquad
\vartheta(F_i):=E_i\;, \quad i=1,2 \;,
\end{equation}
and satisfying $\vartheta^2=id$.

We denote (for obvious reasons) by $U_q(su(2))$ the sub Hopf $*$-algebra
of $U_q(su(3))$ generated by the elements $\{K_1,K_1^{-1},E_1,F_1\}$ and by
$U_q(u(2))$ the Hopf $*$-algebra generated by $U_q(su(2))$,
$K_1K_2^2$ and $(K_1K_2^2)^{-1}$. Notice that   
$K_1K_2^2$ commutes with all elements of $U_q(su(2))$.

\bigskip
Irreducible finite dimensional $*$-representations of $U_q(su(3))$ are classified
by two non-negative integers $n_1,n_2$ (see e.g.~\cite{KS97}). The
representation space $V_{(n_1,n_2)}$ has dimension 
$$
\dim V_{(n_1,n_2)}=\tfrac{1}{2}(n_1+1)(n_2+1)(n_1+n_2+2) \;,
$$
and orthonormal basis
$\ket{n_1,n_2,j_1,j_2,m}$, with labels satisfying the constraints
\begin{equation}\label{eq:constr}
j_i=0,1,2,\ldots,n_i \;,\qquad
\oh (j_1+j_2)-|m|\in\N\;.
\end{equation}
The generators of $U_q(su(3))$ act on $V_{(n_1,n_2)}$ as follows
\begin{align*}
K_1\ket{n_1,n_2,j_1,j_2,m} &:=q^m\ket{n_1,n_2,j_1,j_2,m} \;,\\
K_2\ket{n_1,n_2,j_1,j_2,m} &:=q^{\frac{3}{4}(j_1-j_2)+\frac{1}{2}(n_2-n_1-m)}\ket{n_1,n_2,j_1,j_2,m} \;,\\
E_1\ket{n_1,n_2,j_1,j_2,m} &:=\sqrt{[\tfrac{1}{2} (j_1+j_2)-m][\tfrac{1}{2} (j_1+j_2)+m+1]}
   \,\ket{n_1,n_2,j_1,j_2,m+1} \;,\\
E_2\ket{n_1,n_2,j_1,j_2,m} &:=\sqrt{[\tfrac{1}{2} (j_1+j_2)-m+1]}\,A_{j_1,j_2}\ket{n_1,n_2,j_1+1,j_2,m-\oh} \\
     &\qquad\qquad\qquad +\sqrt{[\tfrac{1}{2} (j_1+j_2)+m]}\,B_{j_1,j_2}\ket{n_1,n_2,j_1,j_2-1,m-\oh}\;,
\end{align*}
with coefficients given by
\begin{align*}
A_{j_1,j_2}&:=\sqrt{\frac{[n_1-j_1][n_2+j_1+2][j_1+1]}
              {[j_1+j_2+1][j_1+j_2+2]}} \;,\\[5pt]
B_{j_1,j_2}&:=\begin{cases}
\sqrt{\dfrac{[n_1+j_2+1][n_2-j_2+1][j_2]}{[j_1+j_2][j_1+j_2+1]}}
\quad & \mathrm{if}\;j_1+j_2\neq 0\;, \\
1 & \mathrm{if}\;j_1+j_2=0\;.
\end{cases}
\end{align*}
As usual, $[z]:=(q^z-q^{-z})/(q-q^{-1})$ denotes the $q$-analogue of $z\in\C$.
The highest weight vector of $V_{(n_1,n_2)}$ is $\ket{n_1,n_2,n_1,0,\oh n_1}$,
corresponding to the weight $(q^{n_1/2},q^{n_2/2})$. There are additional
$*$-representations of  $U_q(su(3))$ that we do not need in the present paper. 
Up to a relabeling, the basis we use is the Gelfand-Tsetlin basis
(Sect.~7.3.3 of~\cite{KS97}). One can pass to the notations of~\cite{Arn97} with
the replacement $E_i=e_i$, $F_i=f_i$, $K_i=q^{h_i/2}$ and
\begin{align*}
n_1&=p_{13}-p_{23}-1 \;,  &
n_2&=p_{23}-p_{33}-1 \;,\\
j_1&=p_{12}-p_{23}-1 \;,  &
j_2&=p_{23}-p_{22} \;, &  2m=2p_{11}-p_{12}-&p_{22}-1 \;.
\end{align*}


The fundamental representation $V_{(0,1)}$ will be needed later on in Sect.~\ref{sec:3}
to construct a pairing of $U_q(su(3))$ with a dual Hopf algebra. Its matrix form,
$\sigma:U_q(su(3))\to\mathrm{Mat}_3(\C)$, is  
\begin{subequations}\label{eq:first}
\begin{align}
\sigma(K_1)&=\ma{
   q^{-1/2} & 0 & 0 \\
   0 & q^{1/2} & 0 \\
   0 & 0 & 1 } \;,&\hspace{-1cm}
\sigma(K_2)&=\ma{
   1 & 0 & 0 \\
   0 & q^{-1/2} & 0 \\
   0 & 0 & q^{1/2} } \;,\\
\sigma(E_1)&=\ma{
   0 & 0 & 0 \\
   1 & 0 & 0 \\
   0 & 0 & 0 } \;,&\hspace{-1cm}
\sigma(E_2)&=\ma{
   0 & 0 & 0 \\
   0 & 0 & 0 \\
   0 & 1 & 0 } \;,
\end{align}
\end{subequations}
having identified  $\ket{0,1,-\oh}$ with $(1,0,0)^t$,
$\ket{0,1,\oh}$ with $(0,1,0)^t$ and $\ket{0,0,0}$ with $(0,0,1)^t$.

\bigskip
In order to have a Casimir operator 
for the algebra $U_q(su(3))$
one needs to enlarge it.
The minimal extension is obtained by adding the element $H:=(K_1K_2^{-1})^{2/3}$
and its inverse; by a slight abuse of notation we continue to use the symbol
$U_q(su(3))$ for this extension.
Such a Casimir element appeared already in \cite[eq. 48]{RP91} but in the framework of formal power series. 
In our notations it reads
\begin{multline}\label{eq:Cq}
\mathcal{C}_q = (q-q^{-1})^{-2} \Big((H+H^{-1})\bigl\{(qK_1K_2)^2+(qK_1K_2)^{-2}\bigr\}
 +H^2+H^{-2}-6 \Big) \\
+(qHK_2^2+q^{-1}H^{-1}K_2^{-2})F_1E_1+(qH^{-1}K_1^2+q^{-1}HK_1^{-2})F_2E_2 \\
+qH[F_2,F_1]_q[E_1,E_2]_q+qH^{-1}[F_1,F_2]_q[E_2,E_1]_q \; ,
\end{multline}
satisfies $\mathcal{C}_q^*=\vartheta(\mathcal{C}_q)=\mathcal{C}_q$ and commutes with all elements of $U_q(su(3))$ as can also be checked by a straightforward computation. Moreover the restriction of $\mathcal{C}_q$ to the irreducible representation $V_{(n_1,n_2)}$ is proportional to the identity (by Schur's lemma) with the constant readily found (by acting on the highest weight vector $v:=\ket{n_1,n_2,n_1,0,\oh n_1}$) to be given by
\begin{equation}\label{SpCq}
\mathcal{C}_q\bigr|_{V_{(n_1,n_2)}}=[\tfrac{1}{3}(n_1-n_2)]^2+
[\tfrac{1}{3}(2n_1+n_2)+1]^2+[\tfrac{1}{3}(n_1+2n_2)+1]^2\;.
\end{equation}

\section{The quantum $SU(3)$ group}\label{sec:3}
The  deformation $\A(SU_q(3))$ of the Hopf $*$-algebra of
representative functions of $SU(3)$ is given in~\cite{RTF90} 
(cf. \cite{KS97}, Sect.~9.2).
As a $*$-algebra it is generated by 9 elements $u^i_j$ ($i,j=1,...,3$)
with commutation relations
\begin{align*}
u^i_ku^j_k &=qu^j_ku^i_k \;,&
u^k_iu^k_j &=qu^k_ju^k_i \;,&&
\forall\;i<j\;, \\
[u^i_l,u^j_k]&=0 \;,&
[u^i_k,u^j_l]&=(q-q^{-1})u^i_lu^j_k \;,&&
\forall\;i<j,\;k<l\;,
\end{align*}
and a cubic relation
$$
\sum\nolimits_{p\in S_3}(-q)^{\ell(\pi)}u^1_{\pi(1)}u^2_{\pi(2)}u^3_{\pi(3)}=1 \;,
$$
where the sum is over all permutations $\pi$ of three elements
and $\ell(\pi)$ is the length of $\pi$.
The $*$-structure is given by
$$
(u^i_j)^*=(-q)^{j-i}(u^{k_1}_{l_1}u^{k_2}_{l_2}-qu^{k_1}_{l_2}u^{k_2}_{l_1})
$$
with $\{k_1,k_2\}=\{1,2,3\}\smallsetminus\{i\}$ and
$\{l_1,l_2\}=\{1,2,3\}\smallsetminus\{j\}$ (as ordered sets).
Thus for example $(u^1_1)^*=u^2_2u^3_3-qu^2_3u^3_2$.
Coproduct, counit and antipode are the usual ones:
$$
\Delta(u^i_j)=\sum\nolimits_ku^i_k\otimes u^k_j\;,\qquad
\epsilon(u^i_j)=\delta^i_j\;,\qquad
S(u^i_j)=(u^j_i)^*\;.
$$

\bigskip
Using the fundamental representation $\sigma:U_q(su(3))\to\mathrm{Mat}_3(\C)$,
given by (\ref{eq:first}), one defines a non-degenerate dual pairing (cf.~\cite{KS97}, Sect.~9.4)
$$
\inner{\,,\,}:U_q(su(3))\times\A(SU_q(3))\to\C\;,\qquad
\inner{h,u^i_j}:=\sigma^i_j(h)\;.
$$
With this pairing -- using Sweedler notation $\Delta(a)=a_{(1)}\otimes a_{(2)}$
for the coproduct --
one gets left and right canonical actions 
$h\az a=a_{(1)}\inner{h,a_{(2)}}$ and $a\za h=\inner{h,a_{(1)}}a_{(2)}$,
explicitly given by
$$
h\az u^i_j=\sum\nolimits_ku^i_k\,\sigma^k_j(h)\;,\qquad
u^i_j\za h=\sum\nolimits_k \sigma^i_k(h)u^k_j\;, 
$$
and which make $\A(SU_q(3))$ an $U_q(su(3))$-bimodule $*$-algebra. It is convenient
to convert the right action into a second left action $\aaz$ commuting  with the action $\az$.
This is done by using the map $\vartheta$ given by (\ref{eq:vartheta}):
$$
h\aaz a:=a\za\vartheta(h) \;,
$$
for all $h\in U_q(su(3))$ and $a\in\A(SU_q(3))$. Since $\vartheta$ is a Hopf
$*$-algebra isomorphism from $U_q(su(3))$ to $U_q(su(3))^{\mathrm{op}}$ the action 
$\aaz$ is compatible with the coproduct and the antipode of  $U_q(su(3))$. Thus, these two left actions make $\A(SU_q(3))$ a left
$U_q(su(3))\otimes U_q(su(3))$-module $*$-algebra.
Explicitly, the actions of generators of $U_q(su(3))$ on generators of $\A(SU_q(3))$ are:
\begin{align*}
K_1\az u^i_1&=q^{-\frac{1}{2}}u^i_1 \;, &
K_1\az u^i_2&=q^{\frac{1}{2}}u^i_2 \;, &
K_1\az u^i_3&=u^i_3 \;, \\
K_2\az u^i_1&=u^i_1 \;, &
K_2\az u^i_2&=q^{-\frac{1}{2}}u^i_2 \;, &
K_2\az u^i_3&=q^{\frac{1}{2}}u^i_3 \;, \\
E_1\az u^i_1&=u^i_2 \;,&
E_1\az u^i_2&=0 \;,&
E_1\az u^i_3&=0 \;,\\
F_1\az u^i_1&=0 \;,&
F_1\az u^i_2&=u^i_1 \;,&
F_1\az u^i_3&=0 \;,\\
E_2\az u^i_1&=0 \;,&
E_2\az u^i_2&=u^i_3 \;,&
E_2\az u^i_3&=0 \;, \\
F_2\az u^i_1&=0 \;,&
F_2\az u^i_2&=0 \;,&
F_2\az u^i_3&=u^i_2 \;,
\end{align*}
and 
\begin{align*}
K_1\aaz u_j^1&=q^{-\frac{1}{2}}u_j^1 \;, &
K_1\aaz u_j^2&=q^{\frac{1}{2}}u_j^2 \;, &
K_1\aaz u_j^3&=u_j^3 \;, \\
K_2\aaz u_j^1&=u_j^1 \;, &
K_2\aaz u_j^2&=q^{-\frac{1}{2}}u_j^2 \;, &
K_2\aaz u_j^3&=q^{\frac{1}{2}}u_j^3 \;, \\
E_1\aaz u_j^1&=u_j^2 \;,&
E_1\aaz u_j^2&=0 \;,&
E_1\aaz u_j^3&=0 \;,\\
F_1\aaz u_j^1&=0 \;,&
F_1\aaz u_j^2&=u_j^1 \;,&
F_1\aaz u_j^3&=0 \;,\\
E_2\aaz u_j^1&=0 \;,&
E_2\aaz u_j^2&=u_j^3 \;,&
E_2\aaz u_j^3&=0 \;, \\
F_2\aaz u_j^1&=0 \;,&
F_2\aaz u_j^2&=0 \;,&
F_2\aaz u_j^3&=u_j^2 \;.
\end{align*}

\bigskip

When computing the spectrum of the `exponential Dirac operator' on $\CP^2_q$ in 
Sect.~\ref{sec:spinc} below, we shall use the fact that the `white' and `black' actions of the Casimir element concide. For the sake of clarity, we state this fact as a Lemma. 
\begin{lemma}\label{eq:azDelta}
Let $\mathcal{C}_q$ be the Casimir element defined in \eqref{eq:Cq}, than 
\begin{equation}
\mathcal{C}_q\az a = \mathcal{C}_q\aaz a\;, \qquad {for ~all} \quad a\in\A(SU_q(3))\,.   
\end{equation}
\end{lemma}
\begin{prova}
Since $\vartheta(\mathcal{C}_q)=\mathcal{C}_q$, this statement is equivalent to 
$\mathcal{C}_q\az a=a \za \mathcal{C}_q$, for all $\;a\in\A(SU_q(3))$, an equality that follows from a simple characterization of the center of $\U$. In fact, if $\U$ and $\A$ are any two Hopf $*$-algebras with a non-degenerate dual pairing $\inner{\,,\,}$ and corresponding left and right canonical actions
$h\az a=a_{(1)}\inner{h,a_{(2)}}$ and $a\za h=\inner{h,a_{(1)}}a_{(2)}$, for $h\in\U$ and $a\in\A$, the center of $\U$ coincides with 
$$
\mathcal{Z}(\U):=\big\{h\in \U\, \big|\,h\az a=a\za h\;,  \forall\;a\in\A \big\}\;.
$$
Indeed, from the definition of the actions,
and non-degeneracy of the pairing, the proposition $\{ h\az a=a\za h \}$ is
equivalent to the proposition
$
\left\{\inner{h'\otimes h,\Delta(a)}=\inner{h\otimes h',\Delta(a)}\;, \forall\;h'\in \U \right\}
$;
this follows from the equalities
$\inner{h'\otimes h,\Delta(a)}=\inner{h',h\az a}$, and $\inner{h\otimes h',\Delta(a)}=\inner{h',a\za h}$.
Then $h\in\mathcal{Z}$ if and only if
$
\inner{h\otimes h',\Delta(a)}=\inner{h'\otimes h,\Delta(a)}$, 
for all $h'\in \U$ and  $a\in\A$. In turn, 
this is equivalent to $\inner{[h,h'],a}=0$, for all $h'\in \U,\; a\in\A$, 
and non-degeneracy of the pairing makes this equivalent to
$[h,h']=0$, for all $h'\in \U$, that is $h$ is in the center of $\U$.
\end{prova}


\bigskip
Below we shall need an explicit basis of `harmonic functions' for the
coordinate algebra on the quantum $5$-sphere, and for some `equivariant
line bundles' on the quantum projective plane.


It follows from general facts  (cf.~Sect.~11 of \cite{KS97}, see also~\cite{Koo94})
that the algebra $\A(SU_q(3))$ is an $U_q(su(3))\otimes U_q(su(3))$-module $*$-algebra and
Peter-Weyl theorem states that it is the multiplicity-free direct
sum of all irreducible representations of $U_q(su(3))\otimes U_q(su(3))$ with highest
weight $(\lambda,\lambda)$, where $\lambda$ runs over all highest
weights of $U_q(su(3))$. These representations are $*$-representations with
respect to the inner product $(a,b)=\varphi(a^*b)$ induced by the
Haar state $\varphi$. `Dually', $\A(SU_q(3))$ is the direct sum of all
its irreducible corepresentations, with multiplicity being the
corresponding dimension. Indeed, we can construct (almost) explicitly the corresponding
`harmonic' orthonormal basis. The element
$\{(u_1^1)^*\}^{n_1}(u_3^3)^{n_2}$ is annihilated by both
$E_i\az$ and $E_i\aaz$ and satisfies $K_i\az\{(u_1^1)^*\}^{n_1}(u_3^3)^{n_2}
=K_i\aaz\{(u_1^1)^*\}^{n_1}(u_3^3)^{n_2}=q^{n_i/2}\{(u_1^1)^*\}^{n_1}(u_3^3)^{n_2}$.
Then the highest weight vector in $\A(SU_q(3))$ corresponding
to the weight $\lambda=(n_1,n_2)$ is
$$
c_{n_1,n_2}\{(u_1^1)^*\}^{n_1}(u_3^3)^{n_2}\;,
$$
with $c_{n_1,n_2}$ a normalization constant.
The remaining vectors of the basis are computed using the following
Lemma.
Recall that the $q$-factorial is defined by
$[n]!:=[n][n-1]\ldots[2][1]$ for $n$ a positive integer, while $[0]!:=1$.
The $q$-binomial is given by
$$
\sqbn{n \\ m}:=\frac{[n]!}{[m]![n-m]!}\;.
$$

\begin{lemma}
With $\ket{n_1,n_2,j_1,j_2,m}$ the basis of the irreducible
representation $V_{(n_1,n_2)}$ of $U_q(su(3))$ described in Sect.~\ref{sec:2}, we have that
$$
\ket{n_1,n_2,j_1,j_2,m}=X_{j_1,j_2,m}^{n_1,n_2}
\ket{n_1,n_2,n_1,0,\oh n_1}\;,
$$
where
$$
X_{j_1,j_2,m}^{n_1,n_2}:=N_{j_1,j_2,m}^{n_1,n_2}\sum_{k=0}^{n_1-j_1}
\frac{q^{-k(j_1+j_2+k+1)}}{[j_1+j_2+k+1]!}\sqbn{n_1-j_1 \\ k}
F_1^{\frac{1}{2}(j_1+j_2)-m+k}[F_2,F_1]_q^{n_1-j_1-k}F_2^{j_2+k}
$$
and
$$
N_{j_1,j_2,m}^{n_1,n_2}:=\sqrt{[j_1+j_2+1]}\sqrt{\frac{[\tfrac{j_1+j_2}{2}+m]!}
{[\tfrac{j_1+j_2}{2}-m]!}\,\frac{[n_2-j_2]![j_1]!}{[n_1-j_1]![j_2]!}\,
\frac{[n_1+j_2+1]![n_2+j_1+1]!}{[n_1]![n_2]![n_1+n_2+1]!}}\;.
$$
\end{lemma}

\begin{prova}
Consider the map $T\in\mathrm{Aut}(V_{(n_1,n_2)})$ defined by
$$
T\ket{n_1,n_2,j_1,j_2,m}=X_{j_1,j_2,m}^{n_1,n_2}\ket{n_1,n_2,n_1,0,\oh n_1} \;.
$$
One checks that $[T,h]v=0$ for any $v\in V_{(n_1,n_2)}$ and
any $h\in U_q(su(3))$. It is enough to do the check for
$h\in\{H_i,E_i,F_i\}_{i=1,2}$. Thus for example, if $h=F_1$ we have
$$
F_1X_{j_1,j_2,m}^{n_1,n_2}=
\sqrt{[\tfrac{j_1+j_2}{2}+m][\tfrac{j_1+j_2}{2}-m+1]}\,X_{j_1,j_2,m-1}^{n_1,n_2}
$$
and
\begin{align*}
F_1T\ket{n_1,n_2,j_1,j_2,m} &=
\sqrt{[\tfrac{j_1+j_2}{2}+m][\tfrac{j_1+j_2}{2}-m+1]}\;X_{j_1,j_2,m-1}^{n_1,n_2}\ket{n_1,n_2,n_1,0,\oh n_1} \\
&=\sqrt{[\tfrac{j_1+j_2}{2}+m][\tfrac{j_1+j_2}{2}-m+1]}\;T\ket{n_1,n_2,j_1,j_2,m-1} \\
&=TF_1\ket{n_1,n_2,j_1,j_2,m} \;.
\end{align*}
The remaining cases are either straightforward (if $h=K_1,K_2$) or can
be derived in a similar manner using the following commutation rules
(proved by induction on $n$):
\begin{align*}
[E_1,F_1^n] &=[n]F_1^{n-1}\, (q-q^{-1})\, (q^{-n+1}K_1^2-q^{n-1}K_1^{-2}) \;, \\
[E_1,[F_2,F_1]_q^n] &=-[n]q^{n-2}[F_2,F_1]_q^{n-1}F_2K_1^{-2} \;, \\
[E_2,F_2^n] &=[n]F_2^{n-1}\, (q-q^{-1})\, (q^{-n+1}K_2^2-q^{n-1}K_2^{-2}) \;, \\
[E_2,[F_2,F_1]_q^n] &=[n]F_1[F_2,F_1]_q^{n-1}K_2^2 \;, \\
F_2F_1^n-q^{-n}F_1^nF_2&=[n]F_1^{n-1}[F_2,F_1]_q \;.
\end{align*}
By Schur's Lemma, $T$ is then proportional to the identity in every irreducible
representation $V_{(n_1,n_2)}$, with some proportionality
constant $A_{n_1,n_2}$.
Since $X_{n_1,0,\frac{1}{2}n_1}^{n_1,n_2}=1$,
$T\ket{n_1,n_2,n_1,0,\oh n_1}=\ket{n_1,n_2,n_1,0,\oh n_1}$
and we deduce that $A_{n_1,n_2}=1$. This means
$$
X_{j_1,j_2,m}^{n_1,n_2}\ket{n_1,n_2,n_1,0,\oh n_1}=
T\ket{n_1,n_2,j_1,j_2,m}=
\ket{n_1,n_2,j_1,j_2,m} \;,
$$
which concludes the proof.
\end{prova}

{}From this Lemma and Peter-Weyl decomposition, we deduce that
an orthonormal basis of $\A(SU_q(3))$ is given by the elements
\begin{equation}\label{eq:GT}
t(n_1,n_2)^{l_1,l_2,k}_{j_1,j_2,m}:=
c_{n_1,n_2}X_{j_1,j_2,m}^{n_1,n_2}\az\Big(X_{l_1,l_2,k}^{n_1,n_2}\aaz
\{(u_1^1)^*\}^{n_1}(u_3^3)^{n_2}\Big)
\end{equation}
and that the linear isometry
\begin{align*}
& \rule{0pt}{16pt}
Q:\A(SU_q(3))\to\bigoplus\nolimits_{(n_1,n_2)\in\N^2}V_{(n_1,n_2)}\otimes V_{(n_1,n_2)}\;, \\
& \rule{0pt}{16pt}
Q(t(n_1,n_2)^{l_1,l_2,k}_{j_1,j_2,m}):=\ket{n_1,n_2,j_1,j_2,m}\otimes\ket{n_1,n_2,l_1,l_2,k} \,
\end{align*}
is an $U_q(su(3))\otimes U_q(su(3))$-module map, that is
for all $h\in U_q(su(3))$
\begin{align*}
Q(h\az t(n_1,n_2)^{l_1,l_2,k}_{j_1,j_2,m})&=
h\ket{n_1,n_2,j_1,j_2,m}\otimes\ket{n_1,n_2,l_1,l_2,k}\;,\\
Q(h\aaz t(n_1,n_2)^{l_1,l_2,k}_{j_1,j_2,m})&=
\ket{n_1,n_2,j_1,j_2,m}\otimes h\ket{n_1,n_2,l_1,l_2,k}\;.
\end{align*}
{}From now on, we will identify $t(n_1,n_2)^{l_1,l_2,k}_{j_1,j_2,m}$
with its image $\ket{n_1,n_2,j_1,j_2,m}\otimes\ket{n_1,n_2,l_1,l_2,k}$.

\section{The quantum projective plane $\CP_q^2$}\label{sec:4}
The quantum complex projective plane, which we denote by $\CP_q^2$, was
studied already in \cite{Mey95} (see also \cite{Wel00}).
The most natural way to come to $\CP^2_q$ is via the $5$-dimensional
sphere $S^5_q$. We shall therefore start by studying the algebra
$\A(S^5_q)$ of coordinate functions on the latter.


The algebra $\A(S^5_q)$ is made of elements of $\A(SU_q(3))$
which are $U_q(su(2))$-invariant,
$$
\A(S^5_q):=\big\{a\in\A(SU_q(3))\,\big|\,h\aaz a=\epsilon(h)a\;\forall\;h\in U_q(su(2))\big\} \;
$$
and, as such, it is the $*$-subalgebra generated by elements $\{u_i^3,\,i=1,\ldots,3\}$ of the last `row'. In~\cite{VS91} it is proved to be
isomorphic, through the identification $z_i=u_i^3$, to the abstract
$*$-algebra with generators $z_i,z_i^*$ and relations:
\begin{gather*}
z_iz_j=qz_jz_i\quad\forall\;i<j \;,\qquad\quad
z_i^*z_j=qz_jz_i^*\quad\forall\;i\neq j \;, \\
\rule{0pt}{16pt}
[z_1^*,z_1]=0 \;,\qquad
[z_2^*,z_2]=(1-q^2)z_1z_1^* \;,\qquad
[z_3^*,z_3]=(1-q^2)(z_1z_1^*+z_2z_2^*) \;, \\
\rule{0pt}{16pt}
z_1z_1^*+z_2z_2^*+z_3z_3^*=1 \;.
\end{gather*}
Since $K_1K_2^2$ is in the commutant of $U_q(su(2))$,
in addition to the `white' action of $U_q(su(3))$,
the algebra $\A(S^5_q)$ carries a `black' action
of the Hopf $*$-algebra generated by $K_1K_2^2$ and its inverse,
which we denote by $U_q(u(1))$.
Thus, $\A(S^5_q)$ is an $U_q(su(3))\otimes U_q(u(1))$-module $*$-algebra.

A vector $\ket{n_1,n_2,l_1,l_2,k}$ of the Gelfand-Tsetlin basis of $V_{(n_1,n_2)}$
is invariant for the action of 
$U_q(su(2))$ if and only if $k=0=(l_1+l_2)/2$. Last equality implies $l_1=l_2=0$.
Then an orthonormal basis of $\A(S^5_q)$ is given by
\begin{equation}\label{eq:t000}
t(n_1,n_2)^{0,0,0}_{j_1,j_2,m}
\end{equation}
where the elements $t$'s are given by (\ref{eq:GT}), with $n_1, n_2$ nonnegative integers and 
labels $j_1, j_2, m$ restricted as in \eqref{eq:constr}. Thus, 
we have the decomposition:
$$
\A(S^5_q)\simeq\bigoplus\nolimits_{(n_1,n_2)\in\N^2}V_{(n_1,n_2)}\;,
$$
into irreducible representations of $U_q(su(3))\otimes U_q(u(1))$, with the generator $K_1K_2^2$ acting on $V_{(n_1,n_2)}$ as $q^{n_2-n_1}$
times the identity map. 



\bigskip
The algebra $\A(\CP^2_q)$ of the quantum projective plane $\CP_q^2$  can be 
defined either as a subalgebra of $\A(S^5_q)$ or (equivalently) as a subalgebra of $\A(SU_q(3))$.
Both versions will be used when constructing (anti)-holomorphic forms on $\CP_q^2$ later on. 
We remind that $K_1K_2^2$ is the generator of the Hopf $*$-algebra denoted $U_q(u(1))$ above. Then, we define
\begin{align*}
\A(\CP^2_q)&:=\big\{a\in\A(S^5_q)\,\big|\,K_1K_2^2\aaz a=a\big\} \\
  &\phantom{:}=\big\{a\in\A(SU_q(3))\,\big|\,h\aaz a=\epsilon(h)a\;,\;\,\forall\;h\in U_q(u(2))\big\} \;.
\end{align*}
The $*$-algebra $\A(\CP^2_q)$ is generated by elements $p_{ij}:=(u_i^3)^*u_j^3=z_i^*z_j$, $j=1,2,3$,
with $*$-structure $(p_{ij})^*=p_{ji}$.
The relations split in commutation rules
\begin{align*}
p_{ii}p_{jk}&=q^{\mathrm{sign}(i-j)+\mathrm{sign}(k-i)}p_{jk}p_{ii}
 &&\mathrm{if}\;i,j,k\;\textrm{are distinct}\;,\\
p_{ii}p_{ij}&=q^{\mathrm{sign}(j-i)+1}p_{ij}
p_{ii}-(1-q^2)\textstyle{\sum_{k<i}}\,q^{6-2k}p_{kk}p_{ij}
 &&\mathrm{if}\;i\neq j\;,\\
p_{ij}p_{ik}&=q^{\mathrm{sign}(k-j)}p_{ik}p_{ij}
 &&\mathrm{if}\;i\notin\{j,k\}\;,\\
p_{ij}p_{jk}&=q^{\mathrm{sign}(i-j)+\mathrm{sign}(k-j)+1}p_{jk}p_{ij}-(1-q^2)\textstyle{\sum_{l<j}}\,p_{il}p_{lk}
 &&\mathrm{if}\;i,j,k\;\textrm{are distinct}\;,\\
p_{ij}p_{ji}&=(1-q^2)\left(\textstyle{\sum_{l<i}}\,p_{jl}p_{lj}
-\textstyle{\sum_{l<j}}\,p_{il}p_{li}\right)
 &&\mathrm{if}\;i\neq j\;,
\end{align*}
(here $\mathrm{sign}(0):=0$)
and `projective plane' conditions
\begin{equation}\label{eq:4.2}
\sum\nolimits_kp_{jk}p_{kl}=p_{jl} \;, \qquad\quad
q^4p_{11}+q^2p_{22}+p_{33}=1 \;.
\end{equation}
The relations above are obtained straightforwardly from those of $\A(S^5_q)$. 
There cannot be additional generators: since $K_1K_2^2\aaz z_i=qz_i$ and
$K_1K_2^2\aaz z_j^*=q^{-1}z_j^*$, a monomial in the algebra $\A(S^5_q)$
is invariant if and only if it contains the same number of $z_i$ and $z_i^*$'s,
which can be reordered using the commutation relations of $\A(S^5_q)$.

The elements $p_{ij}$ are assembled in a $3\times 3$ matrix $P$
which by the first relation in (\ref{eq:4.2}) is an idempotent,
$P^2=P$; it is indeed a projection since $P=P^*$ with the given $*$-structure. By
 the second relation in (\ref{eq:4.2}) it has $q$-trace
$$
\tr_q(P):=q^4p_{11}+q^2p_{22}+p_{33}=1 \,.
$$

At the classical value,  $q=1$, of the parameter, the algebra $\A(\CP^2)$ 
is the algebra of (polynomial) functions on the space of size $3$ rank $1$ complex 
projections. This space is diffeomorphic to the projective plane $\CP^2$ by identifying 
each line in $\C^3$ with the range of a projection.

\noindent
Commutativity of the actions $\az$ and $\aaz$ entails that $\A(\CP^2_q)$ is an
$U_q(su(3))$-module $*$-algebra for the action $\az$ with a decomposition of $\A(\CP^2_q)$ into
irreducible representations of $U_q(su(3))$:
$$
\A(\CP^2_q)\simeq\bigoplus\nolimits_{n\in\N}V_{(n,n)}\;.
$$
Indeed, a vector $t(n_1,n_2)^{0,0,0}_{j_1,j_2,m}\in\A(S^5_q)$ is annihilated by
$K_1K_2^2$ if and only if $n_1=n_2$. Thus an orthonormal basis of $\A(\CP^2_q)$
is given by the vectors:
$$
t(n,n)^{0,0,0}_{j_1,j_2,m} \;,
$$
with $n$ a nonnegative integer and labels $j_1, j_2, m$ again restricted as in \eqref{eq:constr}.


\section{The Dolbeault complex} \label{sec:5}
The algebra inclusion $\A(\CP^2_q)\hookrightarrow\A(S^5_q)$ is a noncommutative
analogue of the $U(1)$-principal bundle $S^5\to\CP^2$,
and as in the classical case, `modules of sections of line bundles over $\CP^2_q$' can be constructed, as equivariant maps, via the characters of $U(1)$. 
For $N\in\Z$, we define
\begin{equation}\label{eq:LNa}
\mathcal{L}_N :=
\big\{a\in\A(S^5_q)\,\big|\,K_1K_2^2\aaz a=q^Na\;\big\} \; ;
\end{equation}
in particular $\mathcal{L}_0=\A(\CP^2_q)$.
Being $\A(S^5_q)$ the subalgebra of $\A(SU_q(3))$ made of
$U_q(su(2))$-invariant elements, $\mathcal{L}_N$ can be equivalently described as
\begin{equation}\label{eq:LNb}
\mathcal{L}_N =
\big\{a\in\A(SU_q(3))\,\big|\,K_1K_2^2\aaz a=q^Na\,,\;
h\aaz a=\epsilon(h)a\,,\;\forall\;h\in U_q(su(2))\big\} \;.
\end{equation}
Each $\mathcal{L}_N$ is a $\A(\CP^2_q)$-bimodule.
Moreover, 
since the actions $\az$ and $\aaz$ commute, each $\mathcal{L}_N$ is a $U_q(su(3))$-equivariant (left) $\A(\CP^2_q)$-module,
that is  it  carries a left action of the crossed product
$\A(\CP^2_q)\rtimes U_q(su(3))$. Using the orthonormal basis 
$\{ t(n_1,n_2)^{0,0,0}_{j_1,j_2,m} \}$ of
$\A(S^5_q)$ given by \eqref{eq:t000} on whose elements the generator $K_1K_2^2$ acts as $q^{n_2-n_1}$ times the identity map, we argue that the vectors
$$
t(n,n+N)^{0,0,0}_{j_1,j_2,m}\quad\mathrm{if}\;N\geq 0\;,
\qquad\mathrm{or}\qquad
t(n-N,n)^{0,0,0}_{j_1,j_2,m}\quad\mathrm{if}\;N< 0\;,
$$
form a linear basis of $\mathcal{L}_N$, with $n\in\N$ and $(j_1,j_2,m)$ satisfying
the usual constraints (\ref{eq:constr}). Thus, we have the decomposition into irreducible representations of $U_q(su(3))$:
\begin{align*}
\mathcal{L}_N &\simeq\bigoplus\nolimits_{n\in\N}V_{(n,n+N)} \;, &&
\hspace{-3cm}\mathrm{if}\;N\geq 0\;,\\
\mathcal{L}_N &\simeq\bigoplus\nolimits_{n\in\N}V_{(n-N,n)} \;, &&
\hspace{-3cm}\mathrm{if}\;N<0\;.
\end{align*}

In the commutative ($q\to 1$) limit,  using the K{\"a}hler structure of $\CP^2$ the Hilbert spaces
of chiral spinors can be written as the completion of $\Omega^{(0,\bullet)}:=\bigoplus_k\Omega^{(0,k)}$,
where $\Omega^{(0,k)}$ are antiholomorphic $k$-forms.
As sections of equivariant vector bundles (see e.g.~Sect.~2.4 of \cite{GS99}),
$\Omega^{(0,0)}$ is isomorphic to $\A(\CP^2)$ and $\Omega^{(0,2)}$ to the
commutative limit of $\mathcal{L}_3$.
Contrary to $0$ and $2$ antiholomorphic forms, $1$-forms
are not associated with the principal bundle
$S^5\to\CP^2$ but rather with the $U(2)$-bundle $SU(3)\to\CP^2$,
via a suitable $2$-dimensional representation of $U(2)$.
For the sake of completeness and clarity, we re-derive these classical results in App.~\ref{sec:A}. 

Motivated by the discussion above, in the deformed, $q\not= 1$, case we \emph{define}
antiholomorphic $0$ and $2$-forms as elements of the bimodules
$$
\Omega^{(0,0)}:=\mathcal{L}_0=\A(\CP^2_q)\;,\qquad
\Omega^{(0,2)}:=\mathcal{L}_3\;.
$$
Instead for $1$-forms we use the $*$-representation
$\tau:U_q(su(2))\to\mathrm{Mat}_2(\C)$ given by
$$
\tau(K_1)=\left(\begin{array}{cc} q^{1/2} & 0 \\ 0 & q^{-1/2} \end{array}\right)\;,\qquad
\tau(E_1)=\left(\begin{array}{cc} 0 & 1 \\ 0 & 0 \end{array}\right)\;,
$$
and define $\Omega^{(0,1)}$ as the equivariant $\A(\CP^2_q)$-bimodule
\begin{multline}\label{hif}
\Omega^{(0,1)}:=\big\{v\in\A(SU_q(3))^2 \, \big| \\
\big| \,
K_1K_2^2\aaz v=q^{\frac{3}{2}}v,\;(h_{(1)}\aaz v)\tau(S(h_{(2)}))
=\epsilon(h)v,\;\forall\;h\in U_q(su(2))\big\} \;.
\end{multline}
That is, $v=(v_+,v_-)\in\A(SU_q(3))^2$ belongs to the subspace
$\Omega^{(0,1)}$ if and only if
\begin{subequations}\label{eq:forms}
\begin{align}
K_1K_2^2\aaz (v_+,v_-) &=q^{\frac{3}{2}}(v_+,v_-) \;, &
K_1\aaz (v_+,v_-) &=(q^{\frac{1}{2}}v_+,q^{-\frac{1}{2}}v_-) \;, \label{eq:fA}\\
\rule{0pt}{14pt}
E_1\aaz (v_+,v_-) &=(0,v_+) \;, &
F_1\aaz (v_+,v_-) &=(v_-,0) \;. \label{eq:fB}
\end{align}
\end{subequations}
Also the bimodule $\Omega^{(0,1)}$ carries a representation of $U_q(su(3))$
given by the `white' action, and its decomposition
into irreducible representations of $U_q(su(3))$ is readily found.
With the basis (\ref{eq:GT}) we find that highest weight vectors
of the spin $1/2$ representation of $U_q(su(2))$ have the form
$t(n_1,n_2)_{j_1,j_2,m}^{1,0,+\frac{1}{2}}$ and
$t(n_1,n_2)_{j_1,j_2,m}^{0,1,+\frac{1}{2}}$ (in the former case
(\ref{eq:constr}) forces $n_1\geq 1$, in the latter $n_2\geq 1$).
These are eigenvectors of $K_1K_2^2\aaz$ with eigenvalue
$q^{n_2-n_1+\frac{3}{2}}$ and $q^{n_2-n_1-\frac{3}{2}}$,
respectively. To get a factor $q^{\frac{3}{2}}$ we need
$n_2=n_1$, resp.~$n_2=n_1+3$. Thus $\Omega^{(0,1)}$
is the linear span of the vectors
$$
\big(
\,t(n,n)_{j_1,j_2,m}^{1,0,+\frac{1}{2}}\,,
\,t(n,n)_{j_1,j_2,m}^{1,0,-\frac{1}{2}}\,
\big)\;,\qquad\big(
\,t(n,n+3)_{j_1,j_2,m}^{0,1,+\frac{1}{2}}\,,
\,t(n,n+3)_{j_1,j_2,m}^{0,1,-\frac{1}{2}}\,
\big)\;,
$$
and we have the decomposition into
irreducible representations of $U_q(su(3))$:
$$
\Omega^{(0,1)}\simeq
\bigoplus\nolimits_{n\geq 1}V_{(n,n)}\oplus
\bigoplus\nolimits_{n\geq 0}V_{(n,n+3)}\;.
$$

We are ready to construct a cochain complex
$$
\Omega^{(0,0)}\stackrel{\deb}{\to}
\Omega^{(0,1)}\stackrel{\deb}{\to}
\Omega^{(0,2)}\to 0 \;,
$$
with dual chain complex
$$
0\leftarrow
\Omega^{(0,0)}\stackrel{\deb^\dag}{\leftarrow}
\Omega^{(0,1)}\stackrel{\deb^\dag}{\leftarrow}
\Omega^{(0,2)} \;.
$$
\begin{prop}\label{prop:cochain}
Let  $X$ and $Y$ be the operators
$$
X:=F_2F_1-2[2]^{-1}F_1F_2\;,\qquad
Y:=E_2E_1-2[2]^{-1}E_1E_2\;.
$$
The maps
\begin{subequations}\label{eq:de}
\begin{align}
\deb\; :\Omega^{(0,0)} &\to\Omega^{(0,1)}     \;, &&
\deb a :=(X^*\aaz a\,,\,E_2\aaz a) \;, \label{eq:deA}\\
\rule{0pt}{14pt}
\deb\; :\Omega^{(0,1)} &\to\Omega^{(0,2)}\;, &&
\deb v :=-E_2\aaz v_+ -Y\aaz v_- \;,\label{eq:deB}\\
\intertext{with $v=(v_+,v_-)$, are well defined
and their composition is $\deb^2=0$, that is
$(\Omega^{(0,\bullet)},\deb)$ is a cochain complex. Similarly, the maps}
\deb^\dag:\Omega^{(0,2)} &\to\Omega^{(0,1)}\;, &&
\deb^\dag b :=(-F_2\aaz b\,,\,-Y^*\aaz b) \;,\label{eq:debarA}\\
\rule{0pt}{14pt}
\deb^\dag:\Omega^{(0,1)} &\to\Omega^{(0,0)}     \;, &&
\deb^\dag v :=X\aaz v_++F_2\aaz v_- \;,\label{eq:debarB}
\end{align}
\end{subequations}
are well defined
and their composition is $(\deb^\dag)^2=0$, that is
$(\Omega^{(0,\bullet)},\deb^\dag)$ is a chain complex.
\end{prop}
\noindent
Before we prove this proposition we remark that Serre relations for $U_q(su(3))$
read 
\begin{equation}\label{eq:SerreB}
E_1Y+X^*E_1=0 \;, \qquad
E_2X^*+YE_2=0 \;.
\end{equation}
Moreover, from the commutation relations of $U_q(su(3))$ we get
\begin{equation}\label{eq:boh}
\begin{array}{rlp{1cm}rl}
K_1K_2^2X^* &= q^{\frac{3}{2}}X^*K_1K_2^2 \;, &&
K_1X^*      &= q^{\frac{1}{2}}X^*K_1 \;, \\
\rule{0pt}{18pt}
K_1K_2^2E_2 &= q^{\frac{3}{2}}E_2K_1K_2^2 \;, &&
K_1E_2      &= q^{-\frac{1}{2}}E_2K_1 \;.
\end{array}
\end{equation}
Later on, we shall also need their coproducts:
\begin{equation}\label{eq:copXY}
\begin{array}{rl}
\Delta X &=X\otimes K_1K_2+(K_1K_2)^{-1}\otimes X+\tfrac{1-q^2}{1+q^2}
          (F_2K_1^{-1}\otimes K_2F_1-K_2^{-1}F_1\otimes F_2K_1) \;, \\
\rule{0pt}{18pt}
\Delta Y &=Y\otimes K_1K_2+(K_1K_2)^{-1}\otimes Y+\tfrac{1-q^2}{1+q^2}
          (E_2K_1^{-1}\otimes K_2E_1-K_2^{-1}E_1\otimes E_2K_1) \;.
\end{array}
\end{equation}

\begin{proof}[Proof of Prop.~\ref{prop:cochain}.]\parindent=0in
We start with $\deb$ and we first prove that it is well defined. For any $a\in
\Omega^{(0,0)}=\A(\CP^2_q)$ we want to show that $(v_+,v_-):=\deb a$
satisfies the defining properties (\ref{eq:forms}) of $\Omega^{(0,1)}$.
Definition (\ref{eq:deA}) gives $v_+=X^*\aaz a$ and $v_-=E_2\aaz a$. These, together with the invariance of $a$, proves that $(v_+,v_-)=\deb a$
satisfies (\ref{eq:fA}). Next, we consider the action of $E_1$ and $F_1$.
As $E_1\aaz a=0$,
$$
v_+=X^*\aaz a=E_1E_2\aaz a=E_1\aaz v_- \;,
$$
and since $[F_1,E_2]=0$ and $F_1\aaz a=0$, we have also
$$
F_1\aaz v_-=E_2F_1\aaz a=0 \;.
$$
Thus, two of conditions (\ref{eq:fB}) are satisfied. {}From relations (\ref{serre}) we get
$$
0=(E_1^2E_2-[2]E_1E_2E_1+E_2E_1^2)\aaz a
 =E_1^2E_2\aaz a=E_1\aaz v_+ \;,
$$
having used $E_1\aaz a=0$. The commutation rule 
$[F_1,X^*]=[2]^{-1}E_2(q^{-1}K_1^2+qK_1^{-2})$ yields:
$$
F_1\aaz v_+=[F_1,X^*]\aaz a=
[2]^{-1}E_2(q^{-1}K_1^2+qK_1^{-2})\aaz a=E_2\aaz a=v_- \;.
$$
Hence, all conditions (\ref{eq:forms}) are proved and the map $\deb$
sends $0$-forms to $1$-forms.

Next, we prove that (\ref{eq:deB}) is well defined,
i.e.~for all $v=(v_+,v_-)$ satisfying (\ref{eq:forms}),
the element $b:=-\deb v=E_2\aaz v_++Y\aaz v_-$ is in $\Omega^{(0,2)}=\mathcal{L}_3$.
It is $U_q(su(2))$-invariant: the first identity in (\ref{eq:SerreB})
gives
$$
E_1\aaz b=E_1E_2\aaz v_++E_1Y\aaz v_-=E_1E_2\aaz v_+-X^*E_1\aaz v_- \;,
$$
and from $E_1\aaz v_-=v_+$ (and using also $E_1\aaz v_+=0$), we get
$$
E_1\aaz b=(E_1E_2-X^*)\aaz v_+=2[2]^{-1}E_2E_1\aaz v_+=0 \;.
$$
Thus $b$ is the highest weight vector of a representation
of $U_q(su(2))$. Using (\ref{eq:fA}), we get 
$$
K_1\aaz b=K_1E_2\aaz v_++K_1Y\aaz v_-=
q^{-\frac{1}{2}}E_2K_1\aaz v_++q^{\frac{1}{2}}YK_1\aaz v_-
=E_2\aaz v_++Y\aaz v_-=b \;.
$$
that is  the  highest weight is zero and $b$ carries the  trivial representation $h\aaz b =\epsilon(h) b$. In a similar fashion one proves that $K_1K_2^2\aaz b=q^3b$. We conclude that  $b\in\mathcal{L}_3$ and (\ref{eq:deB}) maps $1$-forms to
$2$-forms.

To prove that $\deb^2=0$ it is enough to compute the action of
$\deb^2$ on a $0$-form $a$. Composition of (\ref{eq:deA}) and (\ref{eq:deB})
yields
$$
\deb^2a=-(E_2X^*+YE_2)\aaz a \;,
$$
which is zero by (\ref{eq:SerreB}).

We omit the proof for $\deb^\dag$ which goes along similar lines.
\end{proof}

In the commutative case, $\Omega^{(0,\bullet)}$ is a graded associative
graded-commutative algebra.
For $q\neq 1$, we know how to multiply $0$-forms by $1$-forms
and by $2$-forms ($\Omega^{(0,1)}$ and $\Omega^{(0,2)}$
are bimodules for $\A(\CP^2_q)=\Omega^{(0,0)}$), but we still
don't know how to multiply two $1$-forms. Next lemma shows how
to do this.

\begin{lemma}
The product of two $1$-forms $v=(v_+,v_-)$ and $w=(w_+,w_-)$,
defined by
$$
v\wprod w:=\tfrac{2}{[2]}(q^{\frac{1}{2}}v_+w_--q^{-\frac{1}{2}}v_-w_+) \;,
$$
is a $2$-form, that is an element of $\Omega^{(0,2)}=\mathcal{L}_3$.
\end{lemma}
\begin{prova}
Clearly $K_1\aaz (v\wprod w)=v\wprod w$ and $K_2\aaz (v\wprod w)=q^{\frac{3}{2}}v\wprod w$.
Further by (\ref{eq:forms}), which are satisfied by both $v$ and $w$,
\begin{align*}
2^{-1}[2]E_1\aaz(v\wprod w) &=
q^{\frac{1}{2}}(K_1^{-1}\aaz v_+)(E_1\aaz w_-)
+q^{-\frac{1}{2}}(E_1\aaz v_-)(K_1\aaz w_+) \\
\rule{0pt}{14pt} &=v_+w_+-v_+w_+=0 \;.
\end{align*}
Hence $v\wprod w$ is the highest weight vector of the trivial
representation of $U_q(su(2))$, which in particular means that
$F_1\aaz (v\wprod w)=0$.
\end{prova}

If $\omega=(a,v,b)$ is a general element of $\Omega^{(0,\bullet)}$,
with $a$ of degree zero, $v$ of degree $1$ and $b$ of degree $2$,
the algebra structure of $\Omega^{(0,\bullet)}$ is
$$
\omega_1\cdot\omega_2=(a_1a_2\;,\;a_1v_2+v_1a_2\;,\;a_1b_2+b_1a_2+v_1\wprod v_2) \;.
$$
It is easy to see that this product is associative, thus making $\Omega^{(0,\bullet)}$
a graded associative algebra (clearly it is not graded-commutative).
This algebra carries a left action of $U_q(su(3))$: the
white action $\az$ acting on components; it is a module $*$-algebra
for this action. Using the faithful Haar
functional $\varphi$ of $SU_q(3)$
we define a non-degenerate inner product on forms,
\begin{equation}\label{eq:inform}
\inner{\omega_1,\omega_2}:=\varphi(a_1^*a_2)+\varphi(v_{1+}^*v_{2+}
+v_{1-}^*v_{2-})+\varphi(b_1^*b_2) \;,
\end{equation}
with respect to which the action of $U_q(su(3))$ is unitary, that is it
corresponds to a $*$-representation (see Lemma 2.5 of \cite{DDL06}),
and the decomposition $\Omega^{(0,\bullet)}:=\bigoplus_n\Omega^{(0,n)}$ is
orthogonal.

The operators $\deb$ and $\deb^\dag$, being defined via the black action,
clearly commute with the above action of $U_q(su(3))$ on forms.
It also follows from Lemma 2.5 of \cite{DDL06} that $h^*\aaz v=(h\aaz )^\dag v$
for all vectors $v$ with entries in $\A(SU_q(3))$ and with respect to the
inner product coming from the Haar state, 
and this easily implies that $\deb^\dag$ is the Hermitian conjugate of $\deb$.

\begin{prop}\label{prop:bound}
The map $\deb$ is a graded-derivation:
\begin{gather*}
\deb(ab)=a(\deb b)+(\deb a)b \;, \qquad
\deb(av)=(\deb a)\wprod v+a(\deb v) \;, \qquad
\deb(va)=(\deb v)a-v\wprod(\deb a) \;, \\
\intertext{while $\deb^\dag$ satisfy:}
[\deb^\dag,a]v=2[2]^{-1}(F_2\aaz a)v_-+q(X\aaz a)v_+ \;, \qquad
[\deb^\dag,a]c=-q^{\frac{3}{2}}(F_2\aaz a\,,\,F_1F_2\aaz a)c \;,
\end{gather*}
for all $a,b\in\Omega^{(0,0)}$, $v\in\Omega^{(0,1)}$,
$c\in\Omega^{(0,2)}$.
\end{prop}
\begin{prova}
{}From the formula (\ref{eq:copXY}) for the coproducts of $X$ and $Y$,
and by covariance of the action, $w=\deb (ab)$ has components
\begin{align*}
w_+ &=(X^*_{(1)}\aaz a)(X^*_{(2)}\aaz b)
    =(X^*\aaz a)b+a(X^*\aaz b) \;, \\
w_- &=(E_2\aaz a)(K_2\aaz b)+(K_2\aaz a)(E_2\aaz b)
    =(E_2\aaz a)b+a(E_2\aaz b) \;,
\end{align*}
and so $w=(\deb a)b+a(\deb b)$. Next,
\begin{align*}
-\deb(av) &=E_2\aaz(av_+)+Y\aaz(av_-) \\
 &=q^{\frac{1}{2}}(E_2\aaz a)v_++a(E_2\aaz v_+)+
   q^{\frac{1}{2}}(Y\aaz a)v_- +a(Y\aaz v_-)
  +\tfrac{1-q^2}{1+q^2}(E_2\aaz a)(K_2E_1\aaz v_-) \\
 &=-a(\deb v)+2[2]^{-1}q^{-\frac{1}{2}}(E_2\aaz a)v_+
   +q^{\frac{1}{2}}(Y\aaz a)v_- \;; \\
\intertext{but $Y\aaz a=E_2E_1\aaz a=-2[2]^{-1}X^*\aaz a$ and so}
\deb(av) &=a(\deb v)-2[2]^{-1}\bigl\{q^{-\frac{1}{2}}(E_2\aaz a)v_+
   -q^{\frac{1}{2}}(X^*\aaz a)v_-\bigr\} \\ &=a(\deb v)+(\deb a)\wprod v\;.
\end{align*}
Similarly,
\begin{align*}
-\deb(va) &=E_2\aaz(v_+a)+Y\aaz(v_-a) \\
  &=-(\deb v)a+q^{-\frac{1}{2}}v_+(E_2\aaz a)
    +q^{-\frac{1}{2}}v_-(Y\aaz a)-\tfrac{q^{-1}-q}{q^{-1}+q}
    (K_2^{-1}E_1\aaz v_-)(E_2K_1\aaz a) \\
  &=-(\deb v)a+2[2]^{-1}q^{\frac{1}{2}}v_+(E_2\aaz a)+
    q^{-\frac{1}{2}}v_-(Y\aaz a) \\
  &=-(\deb v)a+2[2]^{-1}\bigl\{q^{\frac{1}{2}}v_+(E_2\aaz a)-
    q^{-\frac{1}{2}}v_-(X^*\aaz a)\bigr\} \\
  &=-(\deb v)a+v\wprod(\deb a) \;.
\end{align*}
In the same manner one proves the identities involving
$[\deb^\dag,a]$.
\end{prova}

Hence, the data $(\Omega^{(0,\bullet)},\deb)$ give a left-covariant
differential calculus; it is of `dimension' $2$ since we are considering only
the `antiholomorphic' forms.

\section{The spectral triple}\label{sec:spinc}

One could try to define a  `Dolbeault-Dirac' operator $\D$ on $\CP^2_q$
as $\deb+\deb^\dag$; on a compact K{\"a}hler spin manifold this is proportional to 
the Dirac operator of the Levi-Civita connection. We start with a more general  one,
\begin{equation}\label{eq:DDop}
\D\omega:=(\deb^\dag v,\deb a+s\deb^\dag b,s\deb v) \;,
\end{equation}
where $\omega=(a,v,b)$ is a differential  form, and $s\in\R^+$ is arbitrary for the time being. 
We shall be able to check the 
compact resolvent condition only for $s=\sqrt{[2]/2}$. As shown below, for this 
value the square of the operator $D$ is related to the Casimir $\mathcal{C}_q$ 
of $U_q(su(3))$ given in (\ref{eq:Cq}), and
whose spectrum is in \eqref{SpCq}.

Denote with $\HH_+$ the Hilbert space completion of $\Omega^{(0,0)}\oplus\Omega^{(0,2)}$
and with $\HH_-$ the completion of $\Omega^{(0,1)}$, with respect to the inner
product (\ref{eq:inform}). Let $\HH:=\HH_+\oplus\HH_-$ with grading $\gamma:=1\oplus -1$.

\begin{prop}\label{prop:aimST}
For $s=\sqrt{[2]/2}$ in (\ref{eq:DDop}),
the datum $(\A(\CP^2_q),\HH,\D,\gamma)$ is a $0^+$-dimensional $U_q(su(3))$-equivariant
even spectral triple.
\end{prop}

The aim of this section is to prove this proposition.
We have $\A(\CP^2_q)\subset\A(SU_q(3))$ and $\HH\subset L^2(SU_q(3),\varphi)^4$. 
The diagonal lift of the left regular representation of $\A(SU_q(3))$ to
$L^2(SU_q(3),\varphi)^4$ is bounded, thus 
the representation of $\A(\CP^2_q)$ is bounded too.
By Prop.~\ref{prop:bound},
the commutator $[\D,a]$ acts on forms via left multiplication by elements
of $\A(SU_q(3))$,
$$
[\D,a] \omega=\big([\deb^\dag,a]v,[\deb,a]b+s[\deb^\dag,a]c,s[\deb,a]v\big) \;,
\qquad\forall\;\omega=(b,v,c)\;,
$$
and is bounded for any $a\in\A(\CP^2_q)$.
Equivariance holds because forms are defined as equivariant
$\A(\CP^2_q)$-modules, and the operators $\deb$ and $\deb^\dag$  are  $U_q(su(3))$-invariant.  Last 
step is  to check that $D$ has a compact resolvent: we do this by diagonalizing it, which also guarantees 
the existence a self-adjoint extension.

Classically, the K{\"a}hler Laplacian $\D^2$ is half the
Laplace-Beltrami operator $\Delta=\de\de^\dag+\de^\dag\de$,
which in turn, on a symmetric space 
 is related to the quadratic Casimir of
the symmetry algebra. A similar
property holds in the present case.

\begin{prop}
For $s=\sqrt{[2]/2}$, the operator 
$\Delta_{\deb}:=\D^2$ is given by
$$
\Delta_{\deb}\,\omega=[2]^{-1} (\mathcal{C}_q-2)\aaz\omega \;,
$$
for all $\omega\in\Omega^{(0,\bullet)}$. 
\end{prop}
\begin{prova}
Let $a$ be a $0$-form, $v$ a $1$-form and $b$ a $2$-form.
We need to show that
\begin{subequations}
\begin{align}
 \Delta_{\deb}a& =  \deb^\dag\deb a= [2]^{-1} (\mathcal{C}_q-2)\aaz a \;,\label{eq:provaA} \\
 \Delta_{\deb}b& =  s^2\deb\deb^\dag b= [2]^{-1} (\mathcal{C}_q-2)\aaz b \;,\label{eq:provaB} \\
 \Delta_{\deb}v& = (\deb\deb^\dag+s^2\deb^\dag\deb)v= [2]^{-1} (\mathcal{C}_q-2)\aaz v \;. 
\label{eq:provaC}
\end{align}
\end{subequations}
In the following, when acting with elements of $U_q(su(3))$  on forms, the black action is understood and the symbol $\aaz$ is often omitted. 

\bigskip

\textit{Step 1: proof of (\ref{eq:provaA}).}  {}From the definition (\ref{eq:deA}) and (\ref{eq:debarB}) we have 
$$
\deb^\dag\deb\big|_{\Omega^{(0,0)}}=F_2E_2+XX^* \;,
$$
while, using the invariance of $0$-forms: $K_i\aaz a=a$, $E_1\aaz a = F_1\aaz a =0$, and
neglecting terms that vanish on $0$-forms, we rewrite the restriction of $\mathcal{C}_q$ to
$\Omega^{(0,0)}$ as
$$
\mathcal{C}_q\big|_{\Omega^{(0,0)}}
\simeq 2+[2]F_2E_2+(q[F_2,F_1]_q-[F_1,F_2]_q)E_1E_2
\simeq 2+[2](F_2E_2+XX^*) \;.
$$
This proves (\ref{eq:provaA}).

\bigskip

\textit{Step 2: proof of (\ref{eq:provaB}).}
{}From the definition (\ref{eq:debarA}) and (\ref{eq:deB})
$$
\deb\deb^\dag\big|_{\Omega^{(0,2)}}=E_2F_2+YY^* \;,
$$
while, using the symmetry properties of $2$-forms: $K_1\aaz b=b$,
$K_2\aaz b=q^{\frac{3}{2}}b$ and $E_1\aaz b=F_1\aaz b=0$, and
neglecting terms that vanish on $2$-forms, we have
\begin{align*}
\mathcal{C}_q\big|_{\Omega^{(0,2)}}&\simeq 1+[2]^2+[3]^2+(q^2+q^{-2})F_2E_2
   +([F_2,F_1]_q-q[F_1,F_2]_q)E_1E_2 \\
   &\simeq 1+[2]^2+[3]^2+(q^2+q^{-2})F_2E_2+2^{-1}[2]^2Y^*Y \;.
\end{align*}
To compare the last two equations, we need the commutator
\begin{align*}
[2][Y,Y^*] &=[2]E_2[E_1,Y^*]-2[E_1,Y^*]E_2+[2][E_2,Y^*]E_1-2E_1[E_2,Y^*] \\
&=E_2F_2(qK_1^2+q^{-1}K_1^{-2})-2[2]^{-1}F_2(qK_1^2+q^{-1}K_1^{-2})E_2 \\
&-F_1(K_2^2+K_2^{-2})E_1+2[2]^{-1}E_1F_1(K_2^2+K_2^{-2}) \;, \\
\intertext{which, modulo operators vanishing on $\Omega^{(0,2)}$, becomes}
[2][Y,Y^*] &\simeq [2]E_2F_2-4[2]^{-1}F_2E_2 \;.
\end{align*}
This yields (using $-1+[2]^2+[3]^2=[2]^2[3]$)
\begin{align*}
\Big\{\mathcal{C}_q-\tfrac{[2]^2}{2}\deb\deb^\dag-2\Big\}_{\Omega^{(0,2)}}
   &\simeq -1+[2]^2+[3]^2+(q^2+q^{-2})F_2E_2-2^{-1}[2]^2([Y,Y^*]+E_2F_2) \\
   &\simeq [2]^2\bigl([3]-[E_2,F_2]\bigr)
    \simeq [2]^2\bigl([3]-\tfrac{K_2^2-K_2^{-2}}{q-q^{-1}}\bigr) \simeq 0 \;.
\end{align*}
Then on $2$-forms, $\mathcal{C}_q-2=2^{-1}[2]^2\deb\deb^\dag$, which gives
(\ref{eq:provaB}) if{}f $s^2=[2]/2$.

\bigskip

\textit{Step 3: proof of (\ref{eq:provaC}).}
{}From now on, $s=\sqrt{[2]/2}$ is fixed.
Let $w:=(\deb\deb^\dag+s^2\deb^\dag\deb)v$. Then by definition (\ref{eq:de}),
and using $v_-=F_1\aaz v_+$ and $v_+=E_1\aaz v_-$, we get
\begin{align*}
w_+&=(X^*X+s^2F_2E_2)\aaz v_++(X^*F_2+s^2F_2Y)\aaz v_- \\
\rule{0pt}{16pt}
&=\big\{X^*(X+F_2F_1)+s^2F_2(E_2+YF_1)\big\}\aaz v_+ \;,\\
\rule{0pt}{16pt}
w_-&=(E_2X+s^2Y^*E_2)\aaz v_++(E_2F_2+s^2Y^*Y)\aaz v_- \\
\rule{0pt}{16pt}
&=\big\{E_2(XE_1+F_2)+s^2Y^*(E_2E_1+Y)\big\}\aaz v_-\;.
\end{align*}
Using $K_i\aaz v_+=q^{1/2}v_+$ and $E_1\aaz v_+=0$ we get
$$
[2][X^*,F_2F_1]\aaz v_+=(q^2+q^{-2}-2F_2E_2)\aaz v_+\;,\qquad
[Y,F_1]\aaz v_+=E_2\aaz v_+ \;,
$$
as well as (the action $\aaz v_+$ is omitted)
\begin{align*}
[2][X,X^*] &=[2][X,E_1]E_2-2E_2[X,E_1]+[2]E_1[X,E_2]-2[X,E_2]E_1 \\
  &=F_2(K_1^2+K_1^{-2})E_2-2[2]^{-1}E_2F_2(K_1^2+K_1^{-2}) \\
  &-E_1F_1(qK_2^2+q^{-1}K_2^{-2})+2[2]^{-1}F_1(qK_2^2+q^{-1}K_2^{-2})E_1 \\
  &=2[F_2,E_2]-(q^2+q^{-2})[E_1,F_1]=-2-q^2-q^{-2}=-[2]^2 \;,
\end{align*}
that is $X^*X\aaz v_+=(XX^*+[2])\aaz v_+$. Therefore,
$$
[2]w_+=\big\{(q^2+q^{-2})F_2E_2+[2]XE_1E_2+2[3]\bigr\}\aaz v_+  \;.
$$
In turn, using $K_1\aaz v_-=q^{-1/2}v_-$, $K_2\aaz v_-=qv_-$
and $F_1\aaz v_-=0$ we get
$$
[Y,Y^*]\aaz v_-=v_-\;,\qquad
[X,E_1]\aaz v_-=F_2\aaz v_-\;,\qquad
[Y^*,E_2E_1]\aaz v_-=(1-2[2]^{-1}E_2F_2)\aaz v_- \;.
$$
Thus
$$
[2]w_-=\bigl\{[2]E_2F_2+\bigl\{(q^2+q^{-2})E_2E_1-[2]E_1E_2\bigr\}F_1F_2\big\}\aaz v_- \;.
$$
On the other hand, for the action of $\mathcal{C}_q$ on $v_+$ we get
\begin{align*}
(\mathcal{C}_q-2)\aaz v_+ &=\bigl\{
 2[2]^2-2+(q^2+q^{-2})F_2E_2+(q[F_2,F_1]_q-[F_1,F_2]_q)E_1E_2 \bigr\}\aaz v_+ \\
 &=\bigl\{2[3]+(q^2+q^{-2})F_2E_2+[2]XE_1E_2\bigr\}\aaz v_+=[2]w_+ \;,
\end{align*}
while for the action on $v_-$ we get
$$
(\mathcal{C}_q-2)\aaz v_-=
 \bigl\{2(q^2+q^{-2}-1)+[2]F_2E_2+[F_2,F_1]_q[E_1,E_2]_q+q^2[F_1,F_2]_q[E_2,E_1]_q\bigr\}\aaz v_- \;.
$$
To simplify last equation we need some extra work.
Firstly
\begin{align*}
[[F_2,F_1]_q,[E_1,E_2]_q] &=
 (F_2E_2-q^{-2}E_2F_2)K_1^{-2}-(E_1F_1-q^{-2}F_1E_1)K_2^2 \\
 &\simeq qF_2E_2-q^{-1}E_2F_2+1 \;, \\
[[F_1,F_2]_q,[E_2,E_1]_q] &=
 (F_1E_1-q^{-2}E_1F_1)K_2^{-2}-(E_2F_2-q^{-2}F_2E_2)K_1^2 \\
 &\simeq q^{-2}-q^{-1}E_2F_2+q^{-3}F_2E_2 \;,
\end{align*}
where now the symbol `$\simeq$' means that we are neglecting operators vanishing on $v_-$.
Using these commutation relations, we arrive at
\begin{align*}
(\mathcal{C}_q-2)\aaz v_- &=
 \bigl\{[2]E_2F_2+[E_1,E_2]_q[F_2,F_1]_q+q^2[E_2,E_1]_q[F_1,F_2]_q\bigr\}\aaz v_- \\
 &=\bigl\{[2]E_2F_2+\bigl\{(q^2+q^{-2})E_2E_1-[2]E_1E_2\bigr\}F_1F_2\bigr\}\aaz v_-=[2]w_- \;.
\end{align*}
This concludes the proof.
\end{prova}

{}From now on, $s=\sqrt{[2]/2}$ is fixed.

\begin{lemma}\label{lemma:spec}
The kernel of $\D$ are the constant $0$-forms, while its non-zero
eigenvalues are
\begin{align*}
& \pm\sqrt{\tfrac{2}{[2]}[n][n+2]}
&& \textrm{\textup{with multiplicity}}\;\,(n+1)^3\;,\\
& \pm\sqrt{[n+1][n+2]} \hspace{-1cm}
&& \textrm{\textup{with multiplicity}}\;\,\tfrac{1}{2}n(n+3)(2n+3) \;,
\end{align*}
for all $n\geq 1$.
\end{lemma}

\begin{prova}
The crucial property is $\mathcal{C}_q\aaz\psi=\mathcal{C}_q\az\psi$, cf. Lemma~\ref{eq:azDelta}.
For the action `$\az$' we have a decomposition into irreducible representations of $U_q(su(3))$ as
\begin{align*}
\Omega^{(0,0)} &\simeq\bigoplus\nolimits_{n\geq 0}V_{(n,n)} \;, \\
\Omega^{(0,1)} &\simeq\bigoplus\nolimits_{n\geq 1}V_{(n,n)}\oplus
                      \bigoplus\nolimits_{n\geq 0}V_{(n,n+3)}\;, \\
\Omega^{(0,2)} &\simeq\bigoplus\nolimits_{n\geq 0}V_{(n,n+3)} \;.
\end{align*}
These two observations allow us to compute the spectrum of the operator
$[2]\D^2=(\mathcal{C}_q-2)\aaz$.
Its eigenvalues are $\{0,\alpha_n,\beta_m\}_{n\geq 1,m\geq 0}$, given
with their multiplicities by (cf.~eq.~\eqref{SpCq})
\begin{align*}
& 0 \;,
&& \textrm{mult.}=1 \;,\\
& \alpha_n:=2[n+1]^2-2=2[n][n+2] \;,
&& \textrm{mult.}=2(n+1)^3 \;,\\
& \beta_m:=[m+2]^2+[m+3]^2-1=[2][m+2][m+3]\;,
&& \textrm{mult.}=(m+1)(m+4)(2m+5) \;.
\end{align*}
Since $\D$ is odd, its spectrum is symmetric ($\D v=\lambda v$ implies $\D\gamma v=-\lambda\gamma v$).
Thus, $\ker\D=\ker\D^2$ is the subspace $V_{(0,0)}$ made of constant $0$-forms,
and positive roots and negative roots $\pm\alpha_n^{1/2}$ and
$\pm\beta_m^{1/2}$ appear in the spectrum with the same multiplicity.
\end{prova}

Since the eigenvalues of $\D$ grows exponentially (counting multiplicities),
we conclude that $(\D+i)^{-\epsilon}$ is of trace class for any $\epsilon>0$
and the metric dimension is $0^+$. In particular, $\D$ has compact resolvent.
This concludes the proof of Prop.~\ref{prop:aimST}.

We stress that the spectrum of $\D$ is a $q$-deformation
of the classical one~\cite{GS99}. The Connes' differential calculus
associated with $\D$ is left-covariant, and it would be interesting
to compare it with the first order covariant differential calculi
studied in \cite{Wel00}.

\bigskip
As a byproduct of Lemma \ref{lemma:spec} we compute the 
cohomology $H_{\deb}^\bullet(\CP^2_q)$ of
the complex $(\Omega^{(0,\bullet)},\deb)$ in Prop.~\ref{prop:cochain}. The property that allows us to compute it is an analogue of Hodge decomposition theorem.
We call \emph{harmonic $n$-forms} the collection
$$
\mathfrak{H}^n:=\big\{\omega\in\Omega^{(0,n)}\,\big|\,\D\omega=0\big\} \;. 
$$
Since for a homogeneous form $\omega$,
$\D\omega$ is the sum of two pieces $\deb\omega$ and $\deb^\dag\omega$
with different degree, both must vanish in order for $\D\omega$ to
be zero. Thus, $\omega$ is harmonic if{}f $\deb\omega=\deb^\dag\omega=0$.

\begin{prop}
For all $n$, there is an orthogonal decomposition
\begin{equation}\label{eq:ort}
\Omega^{(0,n)}=\mathfrak{H}^n\oplus\deb\Omega^{(0,n-1)} \oplus\deb^\dag\Omega^{(0,n+1)} \;.
\end{equation}
In particular this means that there is
exactly one harmonic form for each cohomology class:
$$
H_{\deb}^n(\CP^2_q)\simeq\mathfrak{H}^n=\ker\D\big|_{\Omega^{(0,n)}} \;.
$$
\end{prop}
\begin{prova}
Given two forms $\omega_{1}, \omega_{2}$ of degree $n-1$ and $n+1$ respectively, we have that
$$
\inner{\deb\omega_1,\deb^\dag\omega_2}=\inner{\deb^2\omega_1,\omega_2}=0 ,
$$
with the inner product defined in \eqref{eq:inform}. Thus $\deb\Omega^{(0,n-1)}$ and
$\deb^\dag\Omega^{(0,n+1)}$ are orthogonal.

It remains to show that an $n$-form $\eta$ is orthogonal to both
$\deb\Omega^{(0,n-1)}$ and $\deb^\dag\Omega^{(0,n+1)}$ if{}f
it is harmonic. This follows from non-degeneracy of the inner product:
we have
\begin{gather*}
\inner{\eta,\deb\omega_1}=\inner{\deb^\dag\eta,\omega_1}=0 \;, \qquad
\inner{\eta,\deb^\dag\omega_2}=\inner{\deb\eta,\omega_2}=0 \;,
\end{gather*}
for all $\omega_1\in\Omega^{(0,n-1)}$ and $\omega_2\in\Omega^{(0,n+1)}$,
if{}f $\deb\eta=\deb^\dag\eta=0$, that is
if{}f $\eta$ is harmonic. This establishes the orthogonal decomposition in \eqref{eq:ort}.

Forms in the subspace $\mathfrak{H}^n\oplus\deb\Omega^{(0,n-1)}$
are $\deb$-closed by construction. 
On the other hand, 
a $\deb$-closed form $\omega\in\deb^\dag\Omega^{(0,n+1)}$ must be harmonic. Orthogonality of the decomposition forces it to vanish. It follows that
$$
H_{\deb}^n(\CP^2_q)=\big\{\mathfrak{H}^n\oplus\deb\Omega^{(0,n-1)}\big\}
 / \deb\Omega^{(0,n-1)}=\mathfrak{H}^n \;,
$$
and this concludes the proof.
\end{prova}

\noindent
An immediate consequence of this proposition
and of Lemma \ref{lemma:spec}, is that
$$
H_{\deb}^0(\CP^2_q)=\C \;,\qquad
H_{\deb}^1(\CP^2_q)=H_{\deb}^2(\CP^2_q)=0 \;.
$$


\subsection*{Acknowledgements} 
The work of FD was partially supported by the 
`Belgian project IAP - NOSY'. The work of LD and GL was partially supported by the `Italian project PRIN06 - Noncommutative geometry, quantum groups and applications'.
LD acknowledges partial support from the grant PL N201177033. 
\newpage

\appendix
\section{Antiholomorphic forms as equivariant maps}
\label{sec:A}
In this appendix we describe the identification of antiholomorphic forms on $\CP^2$ with suitable equivariant maps on bundles over this manifold. In Sect.~\ref{sec:5} this was  the motivation to define antiholomorphic forms on $\CP^2_q$ as equivariant maps.

We denote by $\A(SU(3))\subset C^\infty(SU(3))$ the algebra of polynomials in the coordinate functions $u=(u^k _j)_{k,j=1,2,3}$ which associate to $g\in SU(3)$ its matrix entries: 
$u^k _j(g):=g^k _j$. Abstractly, $\A(SU(3))$ is the $*$-algebra generated by elements $u^k _j$ for $k,j=1,2,3$, with relations
$$
u^k _j u^h_l=u^h_l u^k _j\;,\qquad
\sum\nolimits_{\pi\in S_3}(-1)^{|\pi|}u^1_{\pi(1)}u^2_{\pi(2)}u^3_{\pi(3)}=1 \;,
$$
where $S_3$ are all permutations $\pi$ of three elements and $|\pi|$ is the sign of $\pi$.
The real structure is 
$$
(u^k _j)^*=(-1)^{j-k}(u^{k_1}_{l_1}u^{k_2}_{l_2}-u^{k_1}_{l_2}u^{k_2}_{l_1}) ,
$$
where $\{k_1,k_2\}=\{1,2,3\}\smallsetminus\{k\}$ and
$\{l_1,l_2\}=\{1,2,3\}\smallsetminus\{j\}$ (as ordered sets).
The above relations are just the statements that $uu^\dag=u^\dag u=1$.
As a Hopf algebra, $\A(SU(3))$ has usual coproduct, counit and antipode,
which are obtained by dualizing the group operations.

The Hopf-algebra $U(sl(3))$ is generated by six elements $H_1,H_2,E_1,E_2,F_1,F_2$
subject to the relations coming from Serre's presentation:
\begin{equation*}
\begin{array}{ccccc}
{}[H_k,E_k]=2E_k\;, &&
{}[H_k,F_k]=-2F_k\;, &&
{}[E_k,F_k]=H_k\;, \\
\rule{0pt}{18pt}
{}[H_k,H_j]=0\;, &&
{}[E_k,F_j]=0\;, &&
{}[H_k,E_j]=-E_j\;, \\
\rule{0pt}{18pt}
{}[H_k,F_j]=F_j\;, &&
(\mathrm{ad}E_k)^2(F_j)=0 \;, &&
(\mathrm{ad}F_k)^2(E_j)=0 \;,
\end{array}
\end{equation*}
for all $k,j=1,2$, with $k\neq j$.
Coproduct, counit and antipode are the standard
ones for a universal enveloping algebra.
The $*$-structure corresponding to the real form $U(su(3))$
of $U(sl(3))$ is given by $H_k^*:=H_k$ and $E_k^*:=F_k$.
The Lie algebra $su(3)$ is recovered as the set of primitive elements satisfying $h^*=-h$. Thus, 
$H_k$, $E_k$ and $F_k$ generate the Lie algebra $sl(3)$ while $su(3)$
is the linear span of $\ii H_k$, $\ii (E_k+F_k)$ and $(E_k-F_k)$.  Of course this is an example of a general statement. 
The Hopf algebra $U(su(3))$ is ({\it cum grano salis}) the `classical limit $q\to 1$' of the Hopf algebra 
$U_q(su(3))$ described in Sect.~\ref{sec:2} and can be obtained from it at the level of formal power series in $\hbar:=\log q$, by setting $K_k=q^{H_k/2}$ and truncating at the $0$-th order in $\hbar$.

The fundamental $*$-representation $\sigma: U(su(3)) \to\mathrm{Mat}_3(\C)$
is given by
\begin{align*}
\sigma(H_1)&= \maa{
   \!\!-1 & 0 & 0 \\
   0 & 1 & 0 \\
   0 & 0 & 0 } \,,&
\sigma(H_2)&=  \maa{
   0 & 0 & 0 \\
   0 & \!\!-1 & 0 \\
   0 & 0 & 1 } \,,&
\sigma(E_1)&=\maa{
   0 & 0 & 0 \\
   1 & 0 & 0 \\
   0 & 0 & 0 } \,,&
\sigma(E_2)&=\maa{
   0 & 0 & 0 \\
   0 & 0 & 0 \\
   0 & 1 & 0 } \,.
\end{align*}
With these, the pairing $\inner{\,,\,}:U(su(3))\times C^\infty(SU(3))\to\C$
defined by
$$
\inner{X,f}=\frac{\de}{\de t}\Big|_{t=0}f(e^{t\sigma(X)})\;,\qquad
\textrm{for all}\;X\in su(3)
$$
becomes $\inner{X,u^k_j}=\sigma^k_j(X)$ on generators, and is non-degenerate
when restricted to $\A(SU(3))$. The actions of $U(su(3))$ on $C^\infty(SU(3))$
via left (resp. right) invariant vector fields are given by
$$
(X\az f)(g)=\frac{\de}{\de t}\Big|_{t=0}f(g\, e^{t\sigma(X)})\;,\qquad
(X\aaz f)(g)=\frac{\de}{\de t}\Big|_{t=0}f(e^{-t\sigma(X)}\, g)\; ,
$$
and are the $q\to 1$ limit of the corresponding actions of $U_q(su(3))$ described in Sect.~\ref{sec:3}, as one can see by computing them for a pair of generators. Note that a left (resp. right) invariant vector field generates a right (resp. left) multiplication on the group but a left (resp. right) action on functions.
In the limit the map $\vartheta$ in (\ref{eq:vartheta}) is
simply the $*$-structure on the real vector space $su(3)$,
extended as a \emph{linear} antimultiplicative map to the whole of  $U(su(3))$; thus
$\sigma(\vartheta(X))=\sigma(X)^*=-\sigma(X)$ for all $X\in su(3)$. 

\bigskip
Functions on the sphere $S^5$ are identified with functions on
$SU(3)$ which are annihilated  by the action $\aaz$ of $H_1,E_1,F_1$. 
They are generated by $z_k =u^3_k , \, k=1,2,3$, for which one has that
$\sum_k z_k z_k ^*=\det(u)=1$. Functions on $\CP^2=S^5/S^1$ are
identified with functions on $S^5$ which are annihilated  by the action $\aaz$ of $H_2$. They are generated by $p_{kj}=z_k ^*z_j$ and correspond to the identification of $\CP^2$ as a real manifold
with the space of $3\times 3$ projections of rank $1$; we
denote $\A(\CP^2)$ the coordinate $*$-algebra generated by $p=(p_{kj})$.

Homogeneous coordinates on $\CP^2$ are classes
$[x_1,x_2,x_3]$, where $(x_1,x_2,x_3)\in\C^3\smallsetminus\{0\}$
and $[x_1,x_2,x_3]=[\lambda x_1,\lambda x_2,\lambda x_3]$
for $\lambda\in\C^*$. We can always choose a representative 
$(x_1,x_2,x_3)\in S^5$.
Local coordinates are given by $x_j/x_k$, in the chart $U_k$ defined by $x_k\neq 0$.
Local coordinate functions $U_k$ are the functions
$\{Z^{(k)}_1,Z^{(k)}_2\}$ associating to each point its local
coordinate, thus
\begin{gather*}
Z^{(1)}_1=z_2/z_1\;,\qquad
Z^{(1)}_2=z_3/z_1\;,\qquad
\mathrm{on}\;\;U_1\;, \\
Z^{(2)}_1=z_1/z_2\;,\qquad
Z^{(2)}_2=z_3/z_2\;,\qquad
\mathrm{on}\;\;U_2\;, \\
Z^{(3)}_1=z_1/z_3\;,\qquad
Z^{(3)}_2=z_2/z_3\;,\qquad
\mathrm{on}\;\;U_3\;,
\end{gather*}
with $z_k $ the generators of $\A(S^5)$. Transition functions are
clearly holomorphic.

An antiholomorphic $1$-form is written as a collection 
$\omega=(\omega^{(j)} )$ and on each chart
$U_j$, 
$$
\omega^{(j)}= f^{(j)}_1\de\bar{Z}^{(j)}_1+f^{(j)}_2\de\bar{Z}^{(j)}_2 \;,
$$
where the coefficients $(f^{(j)}_1,f^{(j)}_2)$ are smooth complex
functions (of $Z^{(j)}_1,Z^{(j)}_2,\bar{Z}^{(j)}_1,\bar{Z}^{(j)}_2$) that must satisfy -- in order for $\omega$ to be uniquely defined -- on each overlap $U_j\cap U_k$, the conditions $(f^{(j)}_1,f^{(j)}_2)g_{jk}=(f^{(k)}_1,f^{(k)}_2)$,
with $g_{jk}:U_j\cap U_k\to GL(2,\C)$ given by
$$
g_{jk}=\ma{
 {\de \bar{Z}_1^{(j)}}/{\de \bar{Z}_1^{(k)}} &
 {\de \bar{Z}_1^{(j)}}/{\de \bar{Z}_2^{(k)}} \\
 {\de \bar{Z}_2^{(j)}}/{\de \bar{Z}_1^{(k)}} &
 {\de \bar{Z}_2^{(j)}}/{\de \bar{Z}_2^{(k)}} } \;.
$$
Explicitly:
$$
g_{12}=g_{21}^{-1}={\bar{z}_2}/{\bar{z}_1^2} \maa{ \!\!-\bar{z}_2 & 0 \\ \!\! -\bar{z}_3 & \bar{z}_1 } \;,\quad
g_{23}=g_{32}^{-1}={\bar{z}_3}/{\bar{z}_2^2} \maa{ \bar{z}_2 & \!\! -\bar{z}_1 \\ 0 & \!\! -\bar{z}_3 } \;,\quad
g_{31}=g_{13}^{-1}={\bar{z}_1}/{\bar{z}_3^2} \maa{ 0 & \!\! -\bar{z}_1 \\ \bar{z}_3 & \!\! -\bar{z}_2 } \;.
$$
The functions $f^{(j)}_k $
can be extended to global $C^\infty$-functions on $\CP^2$ vanishing when $z_j=0$.
For example in the limit $z_2\to 0$, corresponding to $Z_1^{(1)}\to 0$,
the functions $(f^{(1)}_1,f^{(1)}_2)$ are well defined and finite while
$g_{12}$ vanishes; thus from the equality $(f^{(2)}_1,f^{(2)}_2)=(f^{(1)}_1,f^{(1)}_2)g_{12}$
we deduce that $(f^{(2)}_1,f^{(2)}_2)$ vanish too for $z_2\to 0$.
We conclude that, as a $C^\infty(\CP^2)$-bimodule,
\begin{multline*}
\Omega^{(0,1)}\simeq\big\{(f^{(j)}_k )_{i=1,2,\;j=1,2,3}\,\big|\,
f^{(j)}_k \in C^\infty(\CP^2),\;
f^{(j)}_k |_{z_j=0}=0,\; \\
(f^{(j)}_1,f^{(j)}_2)g_{jk}=(f^{(k)}_1,f^{(k)}_2),\;\forall\;i,j,k\big\}.
\end{multline*}

With $\tau$ the spin $1/2$ representation of the algebra $U(su(2))$
generated by $H_1,E_1,F_1$, consider now the $C^\infty(\CP^2)$-bimodule:
\begin{multline*}
\Gamma:=
\big\{v=(v_+,v_-)\in C^\infty(SU(3))^2 \, \big| \\
\big| \,(H_1+2H_2)\aaz v=3v, \;(h_{(1)}\aaz v)\tau(S(h_{(2)}))
=\epsilon(h)v,\;\forall\;h\in U(su(2))\big\} \;;
\end{multline*}
namely, elements of $\Gamma$ are vectors
$v=(v_+,v_-)\in C^\infty(SU(3))^2$ satisfying the conditions:
\begin{subequations}\label{eq:dr}
\begin{align}
H_1\aaz (v_+,v_-)&=(v_+,-v_-) \;,&
E_1\aaz (v_+,v_-)&=(0,v_+) \;,\\
F_1\aaz (v_+,v_-)&=(v_-,0) \;,&
(H_1+2H_2)\aaz (v_+,v_-)&=3(v_+,v_-) \;.
\end{align}
\end{subequations}
The bimodule $\Gamma$ is the $q\to1$ limit of the bimodule in the right hand side of \eqref{hif}. The following result is just the motivation for the identification of that  bimodule as the bimodule of antiholomorphic 1-forms.

\begin{prop}
There is an isomorphism of $C^\infty(\CP^2)$-bimodules
$\psi:\Gamma\to\Omega^{(0,1)}$ given by
$$
(v_+,v_-)\mapsto (f_1^{(j)},f_2^{(j)})=
\psi(v_+,v_-)^{(j)}:=(v_+,v_-)P^{(j)} \;,
$$
where $P^{(j)}\in\mathrm{Mat}_2\bigl(C^\infty(SU(3))\bigr)$ are the following matrices:
$$
P^{(j)}:=\bar{z}_j\ma{
 \;  u^1_k & \;  u^1_l \\
 \! -u^2_k & \! -u^2_l } \;, \qquad j=1,2,3,
$$
with $\{j,k,l\}$ the permutation of $\{1,2,3\}$ with $k<l$.
Under this isomorphism the operator $\deb$ becomes the $q\to 1$ limit
of the operator (\ref{eq:deA}), that is:
\begin{equation}\label{eq:last}
\psi([E_1,E_2]\aaz a,E_2\aaz a)^{(j)}=
\left({\partial a}/{\partial\bar{Z}^{(j)}_k }\,\,,\,
{\partial a}/{\partial\bar{Z}^{(j)}_2}\right) \;,
\end{equation}
for all $a\in C^\infty(\CP^2)$.
\end{prop}
\begin{prova}
Since the algebra of functions is commutative, $\psi$ is a bimodule map (rather than just a left module map).
A priori, $\psi$ maps $(v_+,v_-)\in\Gamma$ into functions $f_k ^{(j)}\in C^\infty(SU(3))$, with $(f_1^{(j)},f_2^{(j)}):=(v_+,v_-)P^{(j)}$. As we shall prove presently, the image of $\psi$ is indeed in $\Omega^{(0,1)}$.

Firstly, the function $f_k ^{(j)}$ vanishes for $z_j=0$, since the
matrix $P^{(j)}$ vanishes there. The relation $(f^{(j)}_1,f^{(j)}_2)g_{jk}=(f^{(k)}_1,f^{(k)}_2)$
follows from the property $P^{(k)}=P^{(j)}g_{jk}$, which is straightforward to check; for instance:
$$
P^{(1)}g_{12}={\bar{z}_2}/{\bar{z}_1}
\ma{
 \!-\bar{z}_2 & \;0 \\
 \! -\bar{z}_3 & \; \bar{z}_1 }
\ma{
 \;  u^1_2 & \;  u^1_3 \\
 \! -u^2_2 & \; -u^2_3 }=
{\bar{z}_2}/{\bar{z}_1}\ma{
 \! -\bar{z}_2u^1_2-\bar{z}_3u^1_3 & \;  \bar{z}_1u^1_3 \\
 \; \bar{z}_2u^2_2+\bar{z}_3u^2_3 & \! -\bar{z}_1u^2_3 } \;,
$$
and the last matrix is just $P^{(2)}$ since $\sum_{k=1}^3\bar{z}_ku^j_k=
\sum_{k=1}^3\bar{u}^3_ku^j_k=0$ for $j\neq 3$.

{}From the properties:
$$
H_1\aaz P^{(j)}=\maa{\!-1 & 0 \\ \;0 & 1}P^{(j)}\;,\quad
E_1\aaz P^{(j)}=\maa{0 & \!-1 \\ 0 & \;0}P^{(j)}\;,\quad
F_1\aaz P^{(j)}=\maa{\;0 & 0 \\ \!-1 & 0}P^{(j)}\;,
$$
and $(H_1+2H_2)\aaz P^{(j)}=-3P^{(j)}$, together with the relations (\ref{eq:dr})
we deduce that the functions $f_k ^{(j)}$ are annihilated by $K_1,K_2,E_1,F_1$: they are functions on $\CP^2$, which proves $\mathrm{Im}(\psi)\subset\Omega^{(0,1)}$.

Since $\det P^{(j)}=(-1)^j(\bar{z}_j)^3$, the matrix $P^{(j)}$
is invertible on $U_j$ and $\psi$ is an isomorphism.

By Leibniz rule, it is enough to prove (\ref{eq:last}) in the case $a=p_{kl}$.
Being $U_1$ dense in $\CP^2$, two forms are equal if{}f they are equal on $U_1$, and  we can also assume $j=1$. On the other hand, $([E_1,E_2],E_2)\aaz p_{kl}=(\bar{u}^1_k,-\bar{u}^2_k)z_l$ and the left
hand side of (\ref{eq:last}) is
$$
\psi([E_1,E_2]\aaz p_{ij},E_2\aaz p_{ij})^{(1)}=(\bar{u}^1_k ,-\bar{u}^2_k )z_jP^{(1)}=
p_{1j}(\delta_{i2}-p_{i2},\delta_{i3}-p_{i3}) \;.
$$
Now, in local coordinates the projection $p=(p_{ij})$ is given by,
$$
p=\frac{1}{1+|Z^{(1)}_1|^2+|Z^{(1)}_2|^2}\ma{1 \\ \bar{Z}^{(1)}_1 \\ \bar{Z}^{(1)}_2}
\Big(1 \;\; Z^{(1)}_1 \;\; Z^{(1)}_2\Big). 
$$
Thus
$$
\frac{\partial p}{\partial\bar{Z}^{(1)}_1}=
{p_{11}}/{p_{21}}\maa{0 \\ & 1 \\ && 0}p-p_{12}p \;,\qquad
\frac{\partial p}{\partial\bar{Z}^{(1)}_2}=
{p_{11}}/{p_{31}}\maa{0 \\ & 0 \\ && 1}p-p_{13}p \;.
$$
To show that
$\,\big(\frac{\partial}{\partial\bar{Z}^{(1)}_1},\frac{\partial}{\partial\bar{Z}^{(1)}_2}\big)p
=\psi([E_1,E_2]\aaz p,E_2\aaz p)^{(1)}\,$ is now a simple algebraic manipulation.
This concludes the proof.
\end{prova}

Since antiholomorphic $2$-forms are the wedge product of two antiholomorphic
$1$-forms, its easy to identify them with equivariant maps: they
are invariant under $U(su(2))$ since $\wedge^2\tau=\epsilon$ is
the trivial representation, while $H_1+2H_2$ acts as multiplication
by $6$, that is: 
$$
\Omega^{(0,2)}\simeq\big\{a\in\A(S^5)\,\big|\,(H_1+2H_2)\aaz a=6a\big\} \;.
$$
This is just the identification of antiholomorphic $2$-forms with the $q\to1$ limit of the bimodule  $\mathcal{L}_3$ defined in \eqref{eq:LNa} as mentioned in Sect.~\ref{sec:5}. 


\end{document}